\documentclass{article}
\usepackage{PRIMEarxiv}

\usepackage[utf8]{inputenc} 
\usepackage[T1]{fontenc}    
\usepackage{hyperref}       
\usepackage{url}            
\usepackage{booktabs}       
\usepackage{amsmath}       
\usepackage{amsfonts}       
\usepackage{nicefrac}       
\usepackage{microtype}      
\usepackage{lipsum}
\usepackage{fancyhdr}       
\usepackage{graphicx}       
\usepackage[figuresright]{rotating}
\usepackage{algorithm}
\usepackage{algpseudocode}

\graphicspath{{media/}}     

\allowdisplaybreaks

\newcommand{\change}[1]{#1}

\title{The Electric Location-Routing Problem: Improved Formulations and Effects of Nonlinear Charging
}

\author{
    Luiz Eduardo Cotta Monteiro \\
    Departamento de Engenharia Industrial \\
    Pontif\'{\i}cia Universidade Cat\'olica do Rio de Janeiro \\
    \texttt{luizmonteiro@aluno.puc-rio.br} \\
    \And
    Rafael Martinelli \\
    Departamento de Engenharia Industrial \\
    Pontif\'{\i}cia Universidade Cat\'olica do Rio de Janeiro \\
    \texttt{martinelli@puc-rio.br} \\
}

\begin{document}
\maketitle

\begin{abstract}
Electric Location-Routing models (ELRP) can contribute to the effective planning of electric vehicles (EVs) fleets and charging infrastructure within EV logistic networks because it simultaneously combines routing and location decisions to find optimal solutions to the network design. This study introduces ELRP models that incorporate nonlinear charging process, multiple charging station types and develop new improved formulations to the problem. Existing ELRP models commonly assume a linear charging process and employ a node-based formulation for tracking EV energy and time consumption. In contrast, we propose novel formulations offering alternative approaches for modeling EV energy, time consumption, and nonlinear charging. Through extensive computational experiments, \change{our analysis demonstrates the effectiveness of the new formulations, reducing the average gap from 29.1\% to 11.9\%, yielding improved solutions for 28 out of 74 instances} compared to the node-based formulation. Moreover, our findings provide valuable insights into the strategic implications of nonlinear charging in ELRP decision-making, offering new perspectives for planning charging infrastructure in EV logistic networks.
\end{abstract}


\section{Introduction}\label{sec:Introduction}

Electric mobility is recognized as a sustainable alternative in the transport sector, as electric vehicles (EVs) can replace internal combustion engine vehicles (ICVs) that use fossil fuels. Recent research indicates an exponential growth trend in the world's EV market accompanied by increased incentive policies for public and private charging infrastructure installations \cite{IEA2024global}. \change{As the production cost of EVs declines and the consumers awareness of low-carbon options increases, vehicle manufacturers tend to increase their investments in EV production while reducing their focus on ICVs \cite{chen2022toward}}.

Despite the sustainability appeal, the development of EV networks still presents many challenges: the cost of batteries is still high, impacting the high price of EVs compared to ICVs \cite{bitencourt2023understanding}; the range covered by EVs is reduced compared to ICVs, and the charging infrastructure is sparse \cite{szumska2021parameters}; the battery recharge time can be excessively long depending on its charging speed mode \cite{gupta2021optimal}; the electric power generation system and distribution grid need to be prepared to support a change in electricity demand patterns related to the EVs \cite{karmaker2024electric}. Due to the limited EV range and the sparse charging infrastructure available, a significant challenge to the expansion of EV logistic networks is related to the effective implementation of charging infrastructure and, in this context, Electric Location-Routing models (ELRP) can be beneficial to find optimized solutions to the problem.

Logistic network design typically combines two optimization problems, the Facility Location Problem and the Vehicle Routing Problem (VRP), to supply a customer network. However, there is an interdependence between facility location and vehicle routing decisions, and treating these decisions independently can lead to sub-optimal solutions \cite{salhi1989effect}. The classical location-routing model (LRP) simultaneously combines location and routing decisions concerning the location of depots and the definition of routes to optimize a logistic network performance \cite{prodhon2014survey}. In the context of EV networks, the ELRP takes another meaning: the charging infrastructure location decision is necessary to keep EVs operational, influencing the viability of routes and the reach of the customer locations. Schiffer and Walther (2018a) name the problem as a Location-Routing problem with intra-route facilities (LRPIF), highlighting the model's decision to locate ``between-route'' facilities \cite{schiffer2018adaptive}. Studies that develop LRP models focusing on electric vehicles are still relatively sparse.

The EV charging infrastructure presented in the ELRP literature has some aspects and approaches to highlight: charging policies differentiate between partial charging, which allows partial battery recharges, and total charging, which forces a full battery recharge \cite{schiffer2017electric}; there may be different types of charging stations \change{(CSs)} which are classified according to charging speed, usually slow, moderate or fast charging; the charging process can be modelled as a nonlinear or linear behavior in the time-energy relationship \cite{montoya2017electric}. Although most ELRP models assume a linear charging function, the actual charging process of EV batteries has a nonlinear behavior due to the variation in voltage and electric current that occurs during the charging process \cite{pelletier2017battery}. All these charging aspects should be considered broadly for a realistic planning of an EV charging infrastructure in an ELRP.

This paper introduces electric location-routing models with nonlinear charging and multiple charging station types (ELRP-NLMS) to support the planning of EV fleets and charging infrastructure in EV logistic networks. The ELRP-NLMS considers simultaneous routing, CS sitting decisions, and the broad aspects of the charging process, including partial charging, multiple CS types, and nonlinear charging. We first extend the classic ELRP formulation proposed by \cite{schiffer2017electric} into a ELRP-NLMS model, and we propose \change{improved formulations} to the problem. The main contributions of this paper are: (1) We introduce the ELRP with the nonlinear charging process and multiple CS types. (2) We develop \change{three improved formulations} for the problem with alternative approaches for modeling EV energy, time consumption, and nonlinear charging processes. (3) We generate \change{instances} for the ELRP-NLMS problem and make them available for future studies. (4) Through extensive computational experiments, our \change{analysis demonstrate the effectiveness of the new formulations, yielding improved solutions for 28 out of 74 instances and reducing the average gap from 29.1\% to 11.9\%} compared to the node-based formulation. (5) Our findings provide valuable insights into the strategic implications of nonlinear charging in ELRP decision-making, offering new perspectives for planning charging infrastructure in EV logistic networks.

The subsequent sections are organized as follows. Section 1.1 provides a comprehensive review of the latest studies pertaining to the ELRP models. The ELRP-NLMS problem definition and its initial model formulation are introduced in Section 2. In Section 3, we propose three novel formulations addressing the problem. Section 4 outlines the computational experiments, the instance generation method, the computational results, and an in-depth discussion. Finally, Section 5 presents the conclusions drawn from the study and avenues for future research.

\subsection{Literature Review}

The first approaches to the interdependence between location and routing decisions refer to the work of \cite{von1961relationship} and \cite{maranzana1964location}. However, the first quantified results that proved the benefits of considering these two decisions simultaneously are credited to \cite{salhi1989effect}. The authors solved seven standard distribution problems from the literature using a two-stage process (location and routing) in a customer network with up to four depots. The results showed that the classical strategy of solving a location problem and a routing problem separately often leads to sub-optimal solutions. The survey research by \cite{prodhon2014survey} indicates a significant increase in the number of papers addressing location-routing problems from 2010 onwards. Comprehensive literature reviews related to LRP models can be found in \cite{schiffer2019vehicle} and \cite{mara2021location}.

To the best of our knowledge, the first problem on modeling simultaneous routing and locating decisions in the context of electric vehicles is the Electric Vehicles Battery Swap Stations Location Routing Problem (BSS-EV-LRP) introduced by \cite{yang2015battery}. The model is formulated as a mixed-integer programming (MIP) model and only considers battery exchange stations. Traditional capacitated vehicle routing instances were adapted, and all nodes are considered potential battery station installation points. \cite{li2015multiple} presents a Multiple charging station location-routing problem with time window, which includes four different types of CSs in the model (fast, moderate, slow and battery swap station). The MIP model solves small instances with six and ten customers, while heuristic methods are applied to solve larger instances. The \change{CSs} accept only the full recharge of the battery.

The Electric location routing problem with time windows and partial recharging (ELRP-TWPR) is proposed by \cite{schiffer2017electric}. Instances are adapted from EVRP models and all customer nodes are considered potential CS installation points. The \change{CSs} allow partial recharges and assumes a linear charging function. The MIP model is solved for instances with five, ten, and fifteen clients and a sensitivity analysis is performed for different objective functions.  \cite{schiffer2018adaptive} develop a generic Location-routing Problem with Intra-route Facilities. The study applies an Adaptative Large Neighborhood Search (ALNS) metaheuristic to solve larger-scale instances in ELRP-TWPR and BSS-EV-LRP problems.
 
\cite{paz2018multi} presents the Multi-depot electric vehicle location routing problem with time windows. In addition to installing two types of CSs, the model also decides on installing multiple depots. One of the CS types is a battery swap station and the CSs assume a linear charging process. The MIP model is solved for instances with up to 15 clients. In the study of \cite{ccalik2021electric}, a Benders decomposition algorithm embedded in a two-phase framework is developed to solve the ELRP model with heterogeneous fleet, formulated as \change{an MIP}. Optimal values are obtained for instances with up to 15 customers. \cite{schiffer2021integrated} develop an integrated framework to evaluate the economic viability of investing in the development of a mid-haul logistic fleet over a multi-period time horizon considering a mixed fleet with EVs and ICVs. The study combines the total cost of ownership analysis with strategic network design decisions in a LRP model. A hybrid metaheuristic with large neighborhood search, local search, and dynamic programming is used to solve the problem. Results show that in certain cases, the electrification of mid-haul logistics is operationally feasible, economically viable and environmentally beneficial.

The Robust electric location-routing problem with time windows and partial recharging (RELRP-TWPR) is introduced by \cite {schiffer2018strategic} as the first ELRP model to consider uncertainty. The RELRP-TWPR is formulated as \change{an MIP} model within an adjustable-robust optimization approach. The uncertainty is related to customer demand, spatial distribution and service time windows. The solution method includes a metaheuristic algorithm and an adversarial approach is applied to solve the robust formulation. \cite{zhang2019novel} proposes the EV battery swap station location-routing problem with stochastic demands model. The model considers only battery swap stations and each customer can randomly assume three types of demand: low, medium or high. Resource policies are applied to deal with uncertainty, in which the EV's route is adjusted to replenishment at the depot or to a battery swap station. A hybrid algorithm is used as a solution method that sequentially applies the Binary Particle Swarm Optimization heuristic for the station location decision and the Variable Neighborhood Search for the routing decision, finding feasible solutions for instances with up to 37 customers.

\cite{guo2022simultaneous} presents the Location-routing model with nonlinear charging and battery wearing. \change{An MIP} model is used to solve ten instances with six to twenty customers; however, none reach the optimal solution. A three-phase algorithm combining the Clarke and Wright algorithm, an iterative greedy algorithm, and an ALNS is proposed to solve larger instances. Despite the study considering the nonlinear charging, the model uses only one type of CS and ignores time tracking along the EV routing and charging process. The study focuses its analysis on the battery wear cost impacts.
\cite{aghalari2023two} presents a Two-stage stochastic location-routing problem for EV fast CS considering stochastic customer demand and ambient temperature fluctuations. The study investigates how the ambient temperature impacts the location-routing decisions of EV fast CS, as temperature variations can influence the EV range. The model considers nonlinear charging but assumes only fixed full recharge for all EVs while ignoring the \change{time constraints and monitoring} along the EV routing and charging process. The MIP model could find feasible solutions for instances with up to 25 customers. \change{A heuristic} Sweep-based Iterated Greedy Adaptative Large Neighborhood algorithm is proposed to solve a case study in Fargo City, North Dakota. The study applies mainly to cities that experience extreme temperature variations. Although the studies by \cite{guo2022simultaneous} and \cite{aghalari2023two} mention the use of nonlinear charging, an ELRP model that includes energy and time aspects of the nonlinear charging process, with partial charging and multiple types of CS has not yet been addressed in the literature.

The literature on exact \change{method} approaches for LRP models is sparse because most studies focused on developing \change{heuristic} methods, as pointed out by \cite{mara2021location}. Among the ELRP models, the exact methods basically focus on developing MIP models and their solutions by commercial solvers. The instance sizes solved by ELRP exact methods are limited, and the largest instance optimally solved has 15 customers, while the largest instance for which a feasible solution was found using exact methods has 30 customers, with a high optimality Gap, above 59\% \cite{schiffer2018strategic}. A summary of the literature related to ELRP formulations and EV charging aspects is described in Table 1, highlighting the contributions of this study. The last columns indicate the mathematical formulation approaches. The Node and Arc tracking columns refer respectively to the formulations that use node and arc-based variables for tracking EV state of chage (SoC) and time consumption, while Recharge Arc refers to the use of a recharge arc variable, including a CS in every EV travel arc; and Path formulation includes multiple \change{CSs} in every EV travel arc. More details about the new ELRP-NLMS formulations and alternative approaches are presented in Section 3.

\begin{table}[htbp]
  \centering
    \caption{Summary of literature related to ELRP models formulation and charging aspects.\label{tab:1}}
    \scriptsize
    \begin{tabular}{p{9.855em}ccccccccccccrccc}
    \toprule
    \multicolumn{1}{c}{} &       & \multicolumn{2}{p{3.43em}}{\textbf{Recharge  Policy}} &       & \multicolumn{3}{p{1.8cm}}{\textbf{Charging Station Type}} &       & \multicolumn{2}{p{1.2cm}}{\textbf{Charging function}} &       & \multicolumn{5}{p{8.495em}}{\textbf{Model Formulation}} \\
\cmidrule{1-1}\cmidrule{3-4}\cmidrule{6-8}\cmidrule{10-11}\cmidrule{13-17}    \textbf{References} &       & \multicolumn{1}{p{2em}}{\begin{sideways}Partial\end{sideways}} & \multicolumn{1}{p{1.43em}}{\begin{sideways}Fixed\end{sideways}} &       & \multicolumn{1}{p{1.57em}}{\begin{sideways}One type\end{sideways}} & \multicolumn{1}{p{2.715em}}{\begin{sideways}Multiple\end{sideways}} & \multicolumn{1}{p{1.855em}}{\begin{sideways}Battery Swap\end{sideways}} &       & \multicolumn{1}{p{1.715em}}{\begin{sideways}Linear\end{sideways}} & \multicolumn{1}{p{1.715em}}{\begin{sideways}Nonlinear  \end{sideways}} &       & \multicolumn{1}{p{1.57em}}{\begin{sideways}Node\end{sideways}} & \multicolumn{2}{p{3.425em}}{\begin{sideways}Arc\end{sideways}} & \multicolumn{1}{p{1.93em}}{\begin{sideways}Recharge Arc\end{sideways}} & \multicolumn{1}{p{1.57em}}{\begin{sideways}Path\end{sideways}} \\
    \multicolumn{1}{c}{} &       &       &       &       &       &       &       &       &       &       &       &       & \multicolumn{1}{p{1.855em}}{Time} & \multicolumn{1}{p{1.57em}}{SoC  } &       &  \\
\cmidrule{1-1}\cmidrule{3-4}\cmidrule{6-8}\cmidrule{10-11}\cmidrule{13-17}    Yang  and Sun (2015) &       &       &       &       &       &       & $\surd$  &       &       &       &       & $\surd$  &       &       &       &  \\
    \midrule
    Li-Ying and Yuan-Bin (2015) &       &       & \multicolumn{1}{l}{$\surd$} &       &       & $\surd$  & $\surd$  &       & $\surd$  &       &       & $\surd$  &       &       &       &  \\
    \midrule
    Schiffer and Walther (2017) &       & $\surd$  &       &       & $\surd$  &       &       &       & $\surd$  &       &       & $\surd$  &       &       &       &  \\
    \midrule
    Paz et al. (2018) &       & $\surd$  &       &       & $\surd$  &       & $\surd$  &       & $\surd$  &       &       & $\surd$  &       &       &       &  \\
    \midrule
    Schiffer and Walther (2018a) &       & $\surd$  &       &       & $\surd$  &       & $\surd$  &       & $\surd$  &       &       & $\surd$  &       &       &       &  \\
    \midrule
    Schiffer and Walther (2018b) &       & $\surd$  &       &       & $\surd$  &       &       &       & $\surd$  &       &       & $\surd$  &       &       &       &  \\
    \midrule
    Schiffer et al. (2021) &       & $\surd$  &       &       &       & $\surd$  &       &       & $\surd$  &       &       & $\surd$  &       &       &       &  \\
    \midrule
    Çalık et al. (2021) &       & $\surd$  &       &       & $\surd$  &       &       &       & $\surd$  &       &       &       &       & $\surd$  &       &  \\
    \midrule
    Guo et al. (2022) &       & $\surd$  &       &       & $\surd$  &       &       &       &       & $\surd$  &       & $\surd$  &       &       &       &  \\
    \midrule
    Aghalari et al. (2023) &       &       & \multicolumn{1}{l}{$\surd$} &       & $\surd$  &       &       &       &       & $\surd$  &       & $\surd$  &       &       &       &  \\
    \midrule
    \textbf{This work} &       & $\surd$  &       &       &       & $\surd$  &       &       & $\surd$  & $\surd$  &       & $\surd$  & $\surd$  & $\surd$  & $\surd$  & $\surd$ \\
    \bottomrule
    \end{tabular}%
  \label{tab:addlabel}%
\end{table}%

\section{The ELRP-NLMS: Problem definition and first model formulation}

Let $C$ be the set of nodes representing the customers, $S$ the set of potential \change{CSs}, $o$ a node representing the depot. The ELRP-NLMS model is defined on a complete and \change{directed} graph $G = (V,A)$ where \change{$V = \{o\} \cup C \cup S$} and $A$ is the set of arcs $(i,j)$ corresponding to road segments connecting nodes $i$ to $j$ where $\{i,j \in V, i \neq j\}$. Each customer $i \in C$ has a service time $sv_i$, and each potential CS $i \in S$ has a charging technology $st_i$ according to slow, moderate, or fast charging speed. The charging process is assumed to be nonlinear and the EV accepts partial and full recharges. The number of \change{CSs} that could be installed is limited to an upper bound $\bar{S}$. An EV traveling from a node $i$ to $j$ implies an energy consumption $e_{ij} \geq 0$ and a driving time $t_{ij} \geq 0$; both the driving time and the energy consumption satisfy the triangular inequality. The customers are serviced by a homogeneous fleet of EVs, while all the EVs have a battery of capacity $Q$, expressed in kWh, and a maximum tour duration $T$. The EVs are assumed to leave the depot with a fully charged battery and the number of EV routes is limited to an upper bound $\bar{R}$. To allow multiple visits to the same CS, a set of dummy nodes $D$ with $\beta$ copies of each CS is used, similar to the formulation presented in \cite{schiffer2017electric}. Each visit of an EV to a CS is modeled as a visit to a distinct copy of the CS, therefore, using CSs copies facilitates tracking the EV battery charge level and the time consumption. The value of $\beta + 1$ corresponds to the number of times each potential CS can be visited by an EV to do a battery recharge. The set $D$ is joined to the set $V$ and we use \change{$S'$} to refer to the subset of \change{CSs} and dummy nodes \change{($S' = S \cup D$)}.

In an ELRP-NLMS feasible solution, each customer must be serviced precisely once by a single EV; each EV route starts and ends at the depot within the maximum duration $T$; each EV route is energy feasible, i.e., the battery level of an EV must remain between 0 and Q throughout the entire route; and \change{the maximum number of CSs that can be installed is equal to $\bar{S}$}. The ELRP-NLMS aims to serve all network customers with a minimum time while respecting all problem constraints.

\subsection{First ELRP-NLMS model formulation}

The first ELRP-NLMS is formulated as a mixed-integer programming model and uses as reference the ELRP model presented by \cite{schiffer2017electric}. However, it incorporates the nonlinear behavior of the charging process as a novelty, using the EVRP model presented by \cite{montoya2017electric} as reference. The MIP uses the following decision variables, $x_{ij}$ is equal to 1 if an EV travels from node $i$ to $j$, and 0 otherwise; $y_{j}$ is equal to 1 if the CS $j$ is opened, and 0 otherwise. Continuous variables $\tau_i$ tracks the EV time, and $q_i$ tracks the EV state of charge (SoC), i.e., the EV battery level, when an EV departs from node $i$. Continuous variables $q^+_j$ and $q^-_j$ specify the SoC when an EV arrives and departs from a CS $j$, \change{while $\Delta_j$ computes the charging time spent in CS $j$}.

The nonlinear charging is modeled by a piecewise linear approximation function using as reference the functions presented by \cite{montoya2017electric}. According to the authors, the functions were adjusted based on data provided by \cite{uhrig2015mobility} and achieved an accurate approximation with a mean absolute error that varies between 0.9\% and 1.9\%. In the ELRP-NLMS models, we developed a new formulation of the charging function that simplifies the modeling, making it more intuitive and with fewer variables and constraints. The piecewise charging function maps the battery levels $q^+_j$ and $q^-_j$ to the equivalent charging times $s_j$ and $d_j$, respectively, to estimate the total charging time $\Delta_j$ spent in CS $j$. This equivalence is computed by a cumulative and proportional sum of the linear segments of the piecewise function.
Let $B$ be the set numbering the piecewise charging function energy segments indexed by $k \in B = \{1, 2\ldots b\}$. The variable $\alpha^+_{jk}$ is the coefficient of each energy segment \change{$E^{'}_{jk}$ for $k \in B$, associated to battery level $q^+_j$, while $\alpha^-_{jk}$ is the coefficient of each energy segment $E^{'}_{jk}$ for $k \in B$, associated to battery level $q^-_j$. The binary variables $z^+_{jk}$ and $z^-_{jk}$ are used to ensure the exact mapping of $q^+_j$ and $q^-_j$ in the charging function. If $q^+_j$ is inside the energy segment $E^{'}_{jk}$, then $\alpha^+_{jk} > 0$, $z^+_{jk} = 1$ and $z^+_{j,k+1} = 0$ for $k \in B\backslash \{b\}$. The same relation is applied to the variables $q^-_j$ , $\alpha^-_{jk}$ and $z^-_{jk}$}. The parameter \change{$T^{'}_{jk}$} refers to the k-th time segment value of the charging function of CS $j$. Figure 1 presents the graphical representation of three charging functions for \change{CSs} with 11, 22 and 44 kW, charging an EV equipped with a 16 kWh battery. \change{Each charging function is composed of three energy and time segments.} The graphic also simulates one EV charging in a moderate speed CS. The first model formulation is described below.

\begin{figure}[htbp]
\caption{Piecewise linear approximation for different types of CS charging (adapted from Montoya et al., 2017).\label{}}
\centerline{\includegraphics[width=13cm]{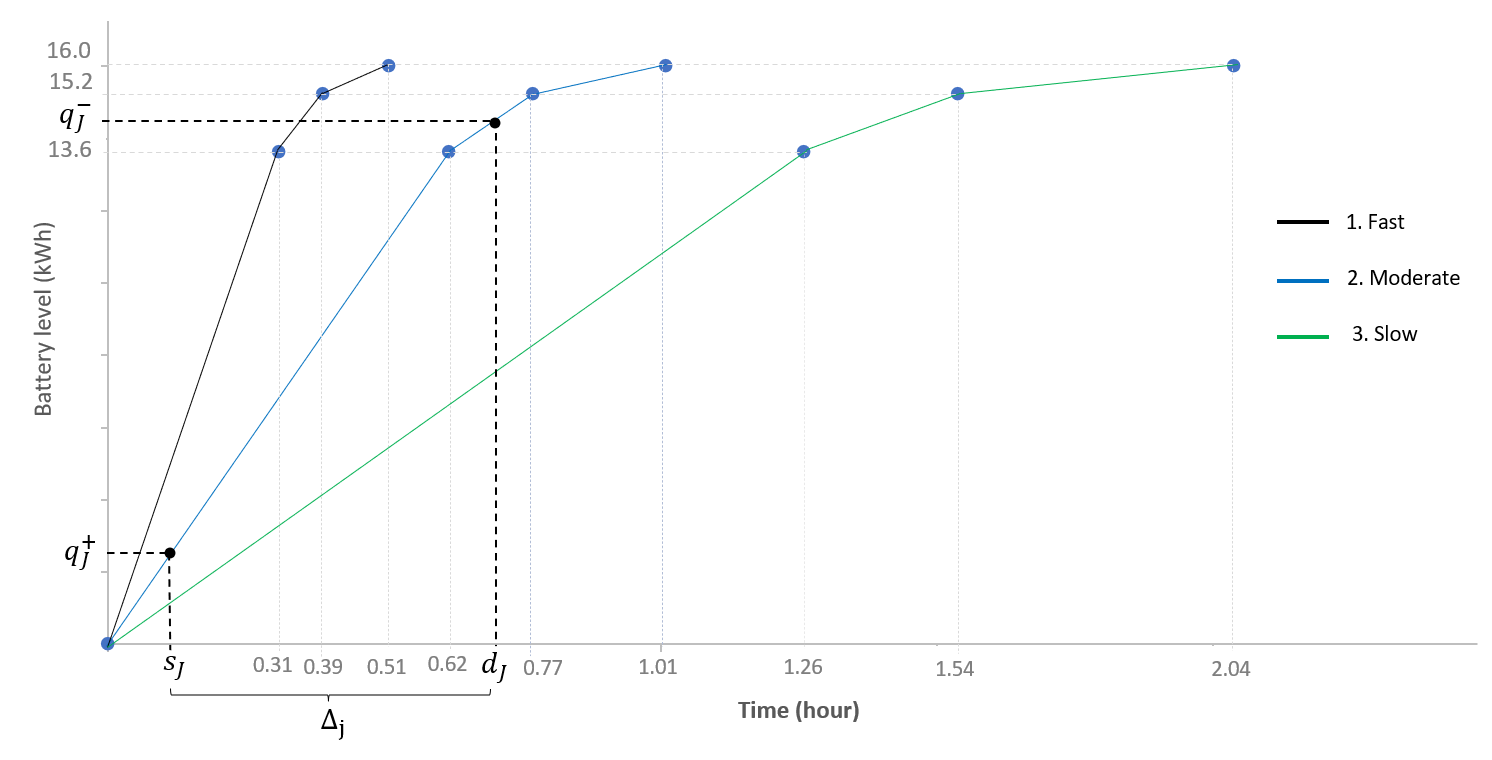}}
\end{figure}

\begin{equation}
    \min \quad \sum\limits_{i \in V} \sum\limits_{j \in V} t_{ij}  x_{ij} + \sum\limits_{j \in \change{S'}} \Delta_j
    \label{eq:m1-of}
\end{equation}

subject to

\begin{align} 
    \sum\limits_{j \in V, i\neq j} x_{ij} &= 1   &     \forall i \in C \label{eq:m1-c1}\\
    \sum\limits_{j \in V, i\neq j} x_{ij} &\leq 1   &     \forall i \in \change{S'} \label{eq:m1-c2}\\
    \sum\limits_{j \in V, i\neq j} x_{ij} &= \sum\limits_{j \in V, i\neq j} x_{ji}   &     \forall i \in V \label{eq:m1-c3}\\
    q_j &\leq q_{i} - e_{ij} x_{ij} + (1 - x_{ij})Q   & \forall i \in V, j \in C, i\neq j \label{eq:m1-c4}\\
    q^+_j &\leq q_{i} - e_{ij} x_{ij} + (1 - x_{ij})Q   &  \forall i \in V, j \in \change{S'}, i\neq j \label{eq:m1-c5}\\
    q_i &\geq e_{io} x_{io}  &  \forall i \in \change{V \backslash \{o\} } \label{eq:m1-c6}\\
    q_o &= Q \label{eq:m1-c7}\\
    q_j &= q^-_j & \forall j \in \change{S'} \label{eq:m1-c8}\\
    q^+_j &= \sum\limits_{k \in B} \alpha^+_{jk} \change{E^{'}}_{jk} & \forall j \in \change{S'} \label{eq:m1-c9}\\
    \alpha^+_{jk+1} &\leq \alpha^+_{jk} & \forall j \in \change{S'}, k \in B \backslash \{b\} \label{eq:m1-c10}\\
    z^+_{jk+1} &\leq \alpha^+_{jk} & \forall j \in \change{S'}, k \in B \backslash \{b\} \label{eq:m1-c11}\\
    z^+_{jk} &\geq \alpha^+_{jk} & \forall j \in \change{S'}, k \in B \label{eq:m1-c12}\\
    s_j &=  \sum\limits_{k \in B} \alpha^+_{jk} \change{T^{'}}_{jk} & \forall j \in \change{S'} \label{eq:m1-c13}\\
    q^-_j &= \sum\limits_{k \in B} \alpha^-_{jk} \change{E^{'}}_{jk} & \forall j \in \change{S'} \label{eq:m1-c14}\\
    \alpha^-_{jk+1} &\leq \alpha^-_{jk} & \forall j \in \change{S'}, k \in B \backslash \{b\} \label{eq:m1-c15}\\
    z^-_{jk+1} &\leq \alpha^-_{jk} & \forall j \in \change{S'}, k \in B \backslash \{b\} \label{eq:m1-c16}\\
    z^-_{jk} &\geq \alpha^-_{jk} & \forall j \in \change{S'}, k \in B \label{eq:m1-c17}\\
    d_j &= \sum\limits_{k \in B} \alpha^-_{jk} \change{T^{'}}_{jk} & \forall j \in \change{S'} \label{eq:m1-c18}\\
    \Delta_j &= d_j - s_j & \forall j \in \change{S'} \label{eq:m1-c19}\\
    \tau_j \geq \tau_i &+ (t_{ij} + sv_{j}) x_{ij} - T (1 - x_{ij}) & \forall i \in V, j \in C \label{eq:m1-c20}\\
    \tau_j \geq \tau_i &+ \Delta_j + t_{ij} x_{ij} - T (1 - x_{ij}) & \forall i \in V, j \in \change{S'} \label{eq:m1-c21}\\
    \tau_j &+ t_{jo}\leq T & \forall j \in V, j \neq o \label{eq:m1-c22}\\
    \tau_{o} &\leq T \label{eq:m1-c23}\\
    q^-_i &- q^+_i \leq Q \, y_i & \forall i \in \change{S'} \label{eq:m1-c24}\\
    \sum\limits_{j \in V} x_{ij} &\leq y_i & \forall i \in \change{S'}, \, i\neq j \label{eq:m1-c25}\\
    \sum\limits_{j \in S} y_j &\leq \bar{S} \label{eq:m1-c26}\\
    \sum\limits_{j \in V} x_{oj} &\leq \bar{R} & \label{eq:m1-c27}\\
    y_i &\geq y_j & \forall \change{{i,j \in S \cup D(S)}, {i < j}} \label{eq:m1-c28}\\
    x_{ij} &= 0 & \forall {i \in \change{S'}}, {j \in \change{S'}, i \neq j} \label{eq:m1-c29}\\
    x_{ij} &\in \{0,1\} & \forall i,j \in V, i \neq j \label{eq:m1-c30}\\
    y_j &\in \{0,1\} & \forall j \in \change{S'} \label{eq:m1-c31}\\
    \tau_i, q_i &\geq 0 & \forall i \in V \label{eq:m1-c32}\\
    q^+_i, q^-_i, s_i, d_i, \Delta_i &\geq 0 & \forall i \in \change{S'} \label{eq:m1-c33}\\
    z^+_{ik}, z^-_{ik} &\in \{0,1\} & \forall i \in \change{S'}, k \in B \label{eq:m1-c34}\\
    0 \leq \alpha^+_{ik}, \alpha^-_{ik} &\leq 1 & \forall i \in \change{S'}, k \in B \label{eq:m1-c35}
\end{align}

The \change{objective function} \eqref{eq:m1-of} aims to minimize the time cost, including EV travel and charging times. \change{This objective function ensure compatibility with the instances presented by \cite{montoya2017electric}, which were extended  to the ELRP-NLMS, and aligns with other ELRP models discussed in Schiffer and Walther (2017) and Paz et al. (2018)}. Constraints \eqref{eq:m1-c1} ensure that each customer is visited once and \eqref{eq:m1-c2} define that each charging station copy is visited at most once. Flow conservation is enforced by \eqref{eq:m1-c3}. Constraints \eqref{eq:m1-c4} and \eqref{eq:m1-c5} track the EV SoC at each customer and each CS, respectively. Constraints \eqref{eq:m1-c6} ensure that when an EV travels to the depot, it must have sufficient energy to complete the route. Constraints \eqref{eq:m1-c7} define the battery charge level when the EV departs from the depot as $Q$, as we assume the EV leaves the depot with a full battery charge. Constraints \eqref{eq:m1-c8} reset the battery tracking to $q^-_j$ when the EV leaves CS $j$.

Constraints \eqref{eq:m1-c9}-\eqref{eq:m1-c13} and \eqref{eq:m1-c14}-\eqref{eq:m1-c18} respectively determine the battery charge level and its corresponding charging time when an EV enters and leaves a CS, based on the piecewise charge function and the respective CS type. Constraints \eqref{eq:m1-c19} compute the time spent at CS $j$. Constraints \eqref{eq:m1-c20} and \eqref{eq:m1-c21} track departure time at each customer and CS, respectively. Constraints \eqref{eq:m1-c22} and \eqref{eq:m1-c23} ensure that any route is completed no later than $T$. CS sitting constraints are defined in \eqref{eq:m1-c24}-\eqref{eq:m1-c26}. Constraints \eqref{eq:m1-c24} ensure that a CS must be built at node $i$ if charging occurs at CS $i$, and \eqref{eq:m1-c25} ensure that only opened \change{CSs} can be part of an EV's route. The number of opened \change{CSs} is limited to $\bar{S}$ by constraints \eqref{eq:m1-c26}. Constraints \eqref{eq:m1-c27} limit the number of EV routes up to $\bar{R}$, defining the size limit of the EV fleet. To avoid symmetry problems related to using dummy nodes as CS copies, Constraint \eqref{eq:m1-c28} defines the order of opened \change{CSs} from the original station to respective dummy nodes \change{and is applied for every node pair $i,j \in S \cup D(S)$, where $D(S)$ represents the set of dummy nodes of each CS in $S$.} Constraint \eqref{eq:m1-c29} removes arcs between dummy vertices of the same CS copy. Constraints \eqref{eq:m1-c30}-\eqref{eq:m1-c35} define the domain of the decision variables.

\subsection{Strengthened first ELRP-NLMS model}

The ELRP combines simultaneous routing and location decisions and is an NP-hard problem as it can be reduced to an EVRP case. Since \change{Model \eqref{eq:m1-of}-\eqref{eq:m1-c35}} becomes computationally difficult, even for small-sized instances, preprocess strategies are used to remove infeasible arcs and strengthen the MIP.
\change{Every EV route should end at the depot or stop at a CS before the battery is empty. The arcs with the minimum energy consumption from each customer to the nearest CS or depot are identified as $\underbar{e}_i$. For nodes $j,i \in C$, the arc $(j,i)$ incident to the customer node $i$ will be infeasible if $\underbar{e}_j + e_{ji} + \underbar{e}_i$ is greater than $Q$, as it exceeds the maximum battery capacity of the electric vehicle. The same logic is applied for the arcs leaving each customer. For $i,j \in C$, the arc (i,j) leaving the customer node $i$ will be infeasible if $\underbar{e}_i + e_{ij} + \underbar{e}_j$ is greater than $Q$.
Infeasible arcs due to battery capacity violations on customer vertices identified in the preprocess are included in the set $I$, therefore $x_{ij}$ is forced to be equal to zero in these conditions using Equations \eqref{eq:m1-c36}. A lower bound for the SoC at each customer is defined by Constraints \eqref{eq:m1-c37}, using the $\underbar{e}_j$ and excluding infeasible values for variables $q_i$ from the model. For the CS nodes $i$, we identify the arcs with the minimum energy consumption ($\underbar{e}^{\prime}_i$) considering the nearest CS or depot and the cyclic tour from $i$ to each customer $j$ ($e_{ij} + e_{ji}$). A lower bound for the final battery charge level at each CS is defined by Constraints \eqref{eq:m1-c38} excluding infeasible values for variables $q^-_i$ and $q_i$ from the model}.

Morevover, to improve the solution of the first ELRP-NLMS model, we apply a procedure based on classic Subtour Elimination Constraints (SEC) by the dynamic \change{introduction} of valid inequalities. The ELRP-NLMS SECs are formulated as Constraints \eqref{eq:m1-c39}, where $A(G^{'})$ is the set of arcs with both end nodes in $G^{'} \subset V$.  For the SEC separation procedure, we refer to \cite{drexl2012note}, \cite{bard2002branch}, and \cite{lysgaard2004new}.

\begin{align}
    x_{ij} &= 0 & \change{\forall (i,j) \in I} \label{eq:m1-c36}\\
        q_i &\geq \underbar{e}_i & \change{\forall i, j \in C} \label{eq:m1-c37}\\
    q^-_i &\geq \underbar{e}^{\prime}_i y_i  & \forall i \in S \label{eq:m1-c38}\\
    \sum\limits_{(i,j) \in A(G^{'})} x_{ij} &\leq |G^{'}| - 1 & \forall G^{'} \subset V \label{eq:m1-c39}
\end{align}

\section{Improved model formulations}

This section presents three new formulations for the ELRP-NLMS. The new formulations are inspired by electric vehicle routing models presented by \cite{froger2019improved} that use arc-based variables and a recharge path concept in the context of EVRP.

\subsection{Second ELRP-NLMS model formulation: Time and energy arc tracking}

The second model formulation substitutes the node-based time and energy tracking variables used in the first model ($\tau_i$ and  $q_i$) by the arc-based continuous variables $\tau_{ij}$ and $q_{ij}$, respectively. Therefore, the modeling strategy for tracking battery and time consumption shifts to a continuous balance across EV travel arcs on routes, where $\tau_{ij}$ represent EV cumulative travel time at node \change{$i$ before} traveling arc $(i,j)$, and $q_{ij}$ the EV battery SoC at node \change{$i$ before} traveling arc $(i,j)$. In the second ELRP-NLMS, constraints (5)-(9), (21)-(24), and (33) of the first model are replaced by the following constraints.

\begin{align} 
    \sum\limits_{j \in V} q_{ij} &= \sum\limits_{l \in V} q_{li} - \change{\sum\limits_{l \in V} e_{li} x_{li} } & & \forall i \in C &  (41)\\
    \sum\limits_{j \in V} q_{ij} &= \change{q^-_i} & & \change{\forall i \in S'}\\
    q^+_j &= \sum\limits_{i \in V} q_{ij} - \change{\sum\limits_{l \in V} e_{ij} x_{ij} } &  & \forall j \in S'\\
    q_{oj} &= \change{Q x_{oj}} & & \forall j \in V \backslash \{o\}\\
    q_{ij} &\leq Q x_{ij} & & \forall i \in V, \, j \in V, \, i \neq j\\
    \tau_{oj} &= \change{0} & & \forall j \in V \backslash \{o\} \\
    \sum\limits_{j \in V} \tau_{ij} &= \change{\sum\limits_{g \in V} \tau_{gi} + \sum\limits_{l \in C} (t_{li} + sv_l) x_{li} + \sum\limits_{l \in S'} (t_{li} x_{li} + \Delta_l)} & & \change{\forall i \in V \backslash \{o\}}\\
    \change{\tau_{io}} &\leq \change{T - t_{io} x_{io} - \Delta_i} & & \change{\forall i \in S'} \\
    \tau_{io} &\leq T \change{- (t_{io} + sv_i) x_{io} } & & \change{\forall i \in C } \\
    \tau_{ij} &\leq T x_{ij} &  & \forall i \in V, \, j \in V, \, i \neq j\\
    \tau_{ij}, \, \; q_{ij} &\geq 0 \;  &  &  \forall i \in V , j \in V, \, i \neq j
\end{align}

Constraints (41) and (42) track the EV battery SoC at node \change{$i$, before} traveling arc $(i,j)$ departing from customer and CS $i$, respectively. Constraints (43) track the battery level when the EV arrives at each CS $j$. Constraints (44) define that the battery charge level when the EV departs from the depot is equal to $Q$, and (45) define that if no EV travels on arc $x_{ij}$, then $q_{ij}$ is 0.

Constraints (46) track EV time \change{when departing from the depot}. Constraints (47) track the EV time \change{before traveling arc $(i,j)$ departing from CS and customer $i$. Constraints (48) and (49) ensure that any route is completed no later than $T$}. Constraints (50) define that if no EV travels on arc $x_{ij}$, then $\tau_{ij}$ is 0. Constraints (51) defines the domain of the $q_{ij}$ and $\tau_{ij}$ variables.

\subsection{Third ELRP-NLMS model formulation: Recharge arc model}

The first two formulations employ dummy nodes to replicate \change{CSs} for multiple recharges. However, the number of CS copies may need to be exceedingly high to guarantee that optimal solutions are not compromised, thereby impacting the size of the model and the time required for solution. In this section, we propose a third formulation for the ELRP-NLMS using recharge arcs approach, which eliminate the need for \change{CSs} replication.

In the recharge arc approach, the variable $x_{ihj}$ is equal to 1 if an EV travels from node $i$ to $j$ performing a recharge at CS $h \in S$, and 0 otherwise. Thus, the charging decision is included in the travel arc variable. To model an EV travel arc that does not pass through a CS, we include a dummy node $d$, where $x_{idj} = 1$ represents an EV travel from $i$ to $j$ without any CS stop. We use $Sd$ to refer to the subset of \change{CSs} and d $(S \cup d)$, and $Co$ to the subset of customers with the depot $(C \cup o)$. The other decision variables are remodeled similarly to $x_{ihj}$ as follows. Continuous variables $\tau_{ihj}$ and $q_{ihj}$ track the EV travel time and SoC of an EV departing from node $i$ traveling to $j$, passing through $h \in Sd$; $q^+_{ihj}$ and $q^-_{ihj}$ track the battery level when an EV arrives and departs from a CS $h$, traveling from node $i$ to $j$; $\Delta_{ihj}$ compute the charging time spent when an EV travels from node $i$ to $j$ passing through the CS $h$. The CS sitting decision variable $y_j$ keeps the same format as the previous models. The third model formulation is as follows.

\begin{equation}
    \min \quad \sum\limits_{i \in Co} \sum\limits_{j \in Co} t_{ij} x_{idj} + \sum\limits_{i \in Co} \sum\limits_{j \in Co} \sum\limits_{h \in S} ( (t_{ih} + t_{hj}) x_{ihj} + \Delta_{ihj})
\end{equation}

subject to

\begin{align} 
    \sum\limits_{j \in Co} \sum\limits_{h \in Sd} x_{ihj} &= 1 & \forall i \in C \\
    \sum\limits_{j \in Co} \sum\limits_{h \in Sd} x_{ihj} &= \sum\limits_{j \in Co} \sum\limits_{h \in Sd} x_{jhi} & \forall i \in Co \\
    q_{ohj} &= Q x_{ohj} & \forall j \in C, h \in Sd\\
    q^+_{ohj} &= Q x_{ohj} - e_{ih} x_{ohj} & \forall j \in C, h \in S\\
    q^+_{ihj} &= q_{ihj} - e_{ih} x_{ihj} & \forall i \in C, j \in C, h \in S\\
    q^+_{ihj} &\leq Q x_{ihj} & \forall i, j \in Co, h \in S\\
    q_{idj} &\leq Q x_{idj} & \forall i, j \in Co\\
    \Delta_{ihj} &\leq \change{\bar{T}_{h}} x_{ihj} & \forall i, j \in Co, h \in S\\
    \begin{split}
        \sum\limits_{l \in Co} \sum\limits_{h \in Sd} q_{ihl} &= \sum\limits_{j \in Co} \sum\limits_{h \in Sd} q_{jhi} - \sum\limits_{j \in Co} e_{ji} x_{jdi}\ -\\ \sum\limits_{j \in Co} &\sum\limits_{h \in S} (e_{jh} + e_{hi}) x_{jhi} + \sum\limits_{j \in Co} \sum\limits_{h \in S} (q^-_{jhi} - q^+_{jhi}) 
    \end{split}&\forall i \in C\\
    q_{iho} &\geq (e_{ih} + e_{ho}) x_{iho} - q^-_{iho} + q^+_{iho} & \forall i \in C, h \in S \\
    q_{ido} &\geq e_{io} x_{ido} & \forall i \in C\\
    q_{ihj} &\leq Q y_{h} & \forall i, j \in Co, h \in S \\
    q^+_{ihj} &= \sum\limits_{k \in B} \alpha^+_{ihjk} \change{E^{'}}_{hk} & \forall i, j \in Co, h \in S\\
    \alpha^+_{ihjk+1} &\leq \alpha^+_{ihjk} & \forall i, j \in Co, h \in S, k \in B \backslash \{b\} \\
    z^+_{ihjk+1} &\leq \alpha^+_{ihjk} & \forall i, j \in Co, h \in S, k \in B \backslash \{b\} \\
    z^+_{ihjk} &\geq \alpha^+{ihjk} & \forall i, j \in Co, h \in S, k \in B\\
    s_{ihj} &=  \sum\limits_{k \in B} \alpha^+_{ihjk} \change{T^{'}}_{hk} & \forall i, j \in Co, h \in S\\
    q^-_{ihj} &= \sum\limits_{k \in B} \alpha^-_{ihjk} \change{E^{'}}_{hk} & \forall i, j \in Co, h \in S\\
    \alpha^-_{ihjk+1} &\leq \alpha^-_{ihjk} & \forall i, j \in Co, h \in S, k \in B \backslash \{b\} \\
    z^-_{ihjk+1} &\leq \alpha^-_{ihjk} & \forall i, j \in Co, h \in S, k \in B \backslash \{b\} \\
    z^-_{ihjk} &\geq \alpha^-_{ihjk} & \forall i, j \in Co, h \in S, k \in B\\
    d_{ihj} &=  \sum\limits_{k \in B} \alpha^-_{ihjk} \change{T^{'}}_{hk} & \forall i, j \in Co, h \in S\\
    \Delta_{ihj} &= d_{ihj} - s_{ihj} & \forall i, j \in Co, h \in S\\
    \begin{split}
    \sum\limits_{l \in Co} \sum\limits_{h \in Sd} \tau_{ihl} &= \sum\limits_{j \in Co} \tau_{jdi} + t_{ji} x_{jdi}\ +\\ \sum\limits_{j \in Co} \sum\limits_{h \in S} &\left(\tau_{jhi} + \Delta_{jhi} + x_{jhi} \left(t_{jh} + t_{hi}\right)\right) + sv_i        
    \end{split} &\forall i \in C\\
    \tau_{ihj} &\leq T x_{ihj} & \forall i, j \in Co, h \in Sd \\
    \tau_{iho} + \Delta_{iho} &\leq (T - t_{ih} - t_{ho}) x_{iho} & \forall i \in C, h \in S \\
    \tau_{ido}  &\leq (T - t_{io}) x_{ido} & \forall i \in C\\
    x_{ihj} &\leq y_h  & \forall i, j \in Co, h \in S \\    
    \sum\limits_{h \in S} y_h &\leq \bar{S} \\
    \sum\limits_{j \in Co} \sum\limits_{h \in Sd} x_{ohj} &\leq \bar{R}\\
    \tau_{ihj}, q_{ihj} &\geq 0 & \forall i, j \in Co, h \in Sd\\
    q^+_{ihj}, q^-_{ihj} &\geq 0 & \forall i, j \in Co, h \in S\\
    s_{ihj}, d_{ihj}, \Delta_{ihj} &\geq 0 & \forall i, j \in Co, h \in S\\
    \alpha^+_{ihjk}, \alpha^-_{ihjk} &\in \left[0, 1\right] & \forall i, j \in Co, h \in S, k \in B\\
    x_{ihj} &\in \{0,1\} & \forall i, j \in Co, h \in Sd\\
    y_h &\in \{0,1\} & \forall h \in S \\
    z^+_{ihjk}, z^-_{ihjk} &\in \{0,1\} & \forall i, j \in Co, h \in S, k \in B
\end{align}

The objective function (52) seeks to minimize the total time. Constraints (53) ensure that each customer is visited once, and (54) enforce the flow conservation. Constraints (55) define the battery SoC when the EV departs from the depot is equal to $Q$. Constraints (56) and (57) track the EV battery SoC when it arrives at CS $h$ after departing from the depot $o$ and from the customer $i$, respectively. Constraints (58) and (59) define that if no EV travels on arc $x_{ihj}$, then $q^+_{ihj}$ and $q_{ihj}$ are equal to 0. Constraints (60) define that if no EV travels on arc $x_{ihj}$, no battery recharge can take on $h \in S$. The parameter \change{$\bar{T}_{h}$} refers to the maximum charging time at CS $h$. Constraints (61) track the battery SoC when EV departs from node $(i \in C)$ to travel arc $x_{ihl}$, and Constraints (62) and (63) ensure that EVs must have sufficient energy to complete the route. Constraints (64) define that $q_{ihj} = 0$ if CS $h$ is not opened.

The formulation of the piecewise linear approximation function follows the same principle as the two previous models, where the continuous variable $\alpha^+_{ihjk}$ is the coefficient of each energy segment $k \in B$ associated to battery level $q^+_{ihj}$ and $\alpha^-_{ihjk}$ is the coefficient of energy segment $k \in B$ associated to battery level $q^-_{ihj}$; the binary variable $z^+_{ihjk}$ is equal to 1 if $q^+_{ihj} \geq \sum_{k \in B} \change{E^{'}}_{ihj(k-1)}$ and $z^-_{ihjk}$ is equal to 1 if $q^-_{ihj} \geq \sum_{k \in B} \change{E^{'}}_{ihj(k-1)}$ where $\change{E^{'}}_{ihjk}$ refers to the $k$-th energy segment value of the piecewise charging function of CS $h$. $\change{T^{'}}_{ihjk}$ refers to the $k$-th time segment value of the charging function of CS $h$. Constraints (65)-(69) and (70)-(74) respectively define the battery charge level and its corresponding charging time when an EV enters and leaves a CS $h$. Constraints (75) compute the time spent at CS $h$.

Constraints (76) track the EV time when it departs from each customer $i$, and (77) define that if no EV travels on arc $x_{ihj}$ then $\tau_{ihj}$ is equal to 0. Constraints (78) and (79) ensure that any route is completed no later than $T$. Constraints (80) ensure that an EV can only travel a recharge arc $x_{ihj}$ for $h \in S$ if the CS $h$ is opened. The number of opened \change{CSs} is limited to $\bar{S}$ by constraints (81), while constraints (82) limit the number of EV routes. Constraints (83)-(89) define the variable domains.

\subsection{Forth ELRP-NLMS model formulation: Recharge path model}

The recharge arc approach restricts up to one CS recharge between two customer nodes, therefore, an optimal solution can be compromised if it includes two consecutive CS recharges. To deal with this limitation, the fourth formulation of the ELRP-NLMS uses the recharge path approach, which considers multiple \change{CSs} between two nodes $(i,j) \in Co$. The formulation includes the set $P_{ij}$ with all possible \change{CS paths} between $(i,j) \in Co$, limited to the upper bound number of installed stations $\bar{S}$. A \change{CS} path $p \in P_{ij}$, is a sequence of different \change{CSs} connecting $i$ to $j$, and a \change{vector $M^p_m = \left[1,\ldots,l\right]$ is the vector of indices where each index $m$ represents the position of a CS along the path $p$. The variables ${q^+}_{ij}^{p,m}$ and ${q^-}_{ij}^{p,m}$ track the battery level when an EV arrives and departs from a CS at position $m$ in path $p$. For the first CS in path $p$, we use the index $f$, so the variable ${q^+}_{ij}^{p,f}$ tracks the battery level when EV arrives at the first CS in $p$. If the vector $M^p_m$ is empty, the EV travels from $i$ to $j$ without any CS stop.}

The variable $x_{ij}^p$ takes the value one if an arc starting from node $i$ traveling path $p$ to node $j$ is traveled by an EV.  The other decision variables are modeled similarly: $\tau_{ij}^p$ and $q_{ij}^p$ track the EV travel time and SoC of an EV departing from node $i$ traveling path $p$ to node $j$. Variable $\Delta_{ij}^{p}$ compute the charging time spent when an EV travels path $p$ from $i$ to $j$. The CS sitting decision variable $y_h$ keeps the same format as the previous models. Moreover, the parameters $t^p_{ij}$ and $e^p_{ij}$ represent respectively the sum of the arcs parameters $e_{ij}$ and $t_{ij}$ within the path $p$. The fourth model formulation is as follows.

\begin{equation}
    \min \quad \sum\limits_{i \in Co} \sum\limits_{j \in Co}  \sum\limits_{p \in P_{ij}} ( t_{ij}^p x_{ij}^p + \Delta_{ij}^p)
\end{equation}

subject to

\begin{align} 
    \sum\limits_{j \in Co} \sum\limits_{p \in P_{ij}} x_{ij}^p &= 1 & \forall i \in C \\
    \sum\limits_{j \in Co} \sum\limits_{p \in P_{ij}} x_{ij}^p &= \sum\limits_{j \in Co} \sum\limits_{p \in P_{ij}} x_{ji}^p & \forall i \in Co\\
    q_{oj}^p &= Q x_{oj}^p & \forall j \in C, p \in P_{ij}\\
    {q^+}_{oj}^{\change{p,f}} &= Q x_{oj}^p - e_{oj}^{\change{p,f}} x_{oj}^p & \forall j \in C, p \in P_{ij}, |p| \geq 1\\
    {q^+}_{ij}^{\change{p,f}} &= q_{ij}^p - e_{ij}^{\change{p,f}} x_{ij}^p & \forall i \in C, j \in \change{Co}, p \in P_{ij}, |p| \geq 1\\
    {q^+}_{ij}^\change{{p,m}} &= {q^-}_{ij}^\change{{p,m-1}} - e_{ij}^\change{{p,m}} x_{ij}^p 
    \begin{split}
        &\forall i, j \in Co, m \in \change{M^p} \backslash \{1\},\\
        &p \in P_{ij}, |p| \geq 2
    \end{split}\\
    \begin{split}
        \sum\limits_{g \in Co} \sum\limits_{p \in P_{ig}} q_{ig}^p &= \sum\limits_{j \in Co} \sum\limits_{p \in P_{ji}} q_{ji}^p - \sum\limits_{j \in Co} \sum\limits_{p \in P_{ji}} e_{ji}^p x_{ji}^p\ +\\ &\sum\limits_{j \in Co} \sum\limits_{m \in \change{M^p}} \sum\limits_{p \in P_{ji}} {q^-}_{ji}^{p,m} - {q^+}_{ji}^\change{{p,m}}
    \end{split} & \forall i \in C\\
    q_{ij}^p &\leq Q x_{ij}^p & \forall i, j \in Co, p \in P_{ij}\\
    {q^+}_{ij}^{p} &\leq Q x_{ij}^p & \forall i, j \in Co, p \in P_{ij}: |p|\geq 1\\
    \Delta_{ij}^p &\leq \change{\bar{T}^{p}} x_{ij}^p & \forall i, j \in Co, p \in P_{ij}: |p|\geq 1\\
    q_{io}^p &\geq x_{io}^p e_{io}^{p} - \sum\limits_{m \in \change{M^p}} ({q^-}_{io}^\change{{p,m}} - {q^+}_{io}^\change{{p,m}}) & \forall i \in C, p \in P_{io}\\
    q_{ij}^p &\leq Q y_{h} & \forall i, j \in Co, p \in P_{ij}, h \in p\\
    {q^+}_{ij}^\change{{p,m}} &= \sum\limits_{k \in B} {\alpha^+}_{ijk}^{\change{p,m}} \change{E^{'p,m}}_{ijk} & \forall p \in P_{ij}, m \in \change{M^p}\\
    {\alpha^+}_{ij(k+1)}^\change{{p,m}} &\leq {\alpha^+}_{ijk}^\change{p,m} & \forall k \in B \backslash \{b\}, p \in P_{ij}, m \in \change{M^p} \\
    {z^+}_{ij(k+1)}^\change{{p,m}} &\leq {\alpha^+}_{ijk}^\change{{p,m}} & \forall k \in B \backslash \{b\}, p \in P_{ij}, m \in \change{M^p} \\
    {z^+}_{ijk}^\change{{p,m}} &\geq {\alpha^+}_{ijk}^\change{{p,m}} & \forall k \in B \backslash \{b\}, p \in P_{ij}, m \in \change{M^p} \\  
    s_{ij}^\change{{p,m}} &= \sum\limits_{k \in B} {\alpha^+}_{ijk}^\change{{p,m}} \change{T^{'p,m}}_{ijk} & \forall p \in P_{ij}, m \in \change{M^p} \\
    {q^-}_{ij}^\change{{p,m}} &= \sum\limits_{k \in B} {\alpha^-}_{ijk}^\change{{p,m}} \change{E^{'p,m}}_{ijk} & \forall p \in P_{ij}, m \in \change{M^p}\\
    {\alpha^-}_{ij(k+1)}^\change{{p,m}} &\leq {\alpha^-}_{ijk}^\change{{p,m}} & \forall k \in B \backslash \{b\}, p \in P_{ij}, m \in \change{M^p} \\
    {z^-}_{ij(k+1)}^\change{{p,m}} &\leq {\alpha^-}_{ijk}^\change{{p,m}} & \forall k \in B \backslash \{b\}, p \in P_{ij}, m \in \change{M^p} \\
    {z^-}_{ijk}^\change{{p,m}} &\geq {\alpha^-}_{ijk}^\change{{p,m}} & \forall k \in B \backslash \{b\}, p \in P_{ij}, m \in \change{M^p} \\  
    d_{ij}^\change{{p,m}} &= \sum\limits_{k \in B} {\alpha^-}_{ijk}^\change{{p,m}} \change{T^{'p,m}}_{ijk} & \forall p \in P_{ij}, m \in \change{M^p}\\
    \Delta_{ij}^{p} &=  \sum\limits_{m \in \change{M^p}} (d_{ij}^{p,m} - s_{ij}^\change{{p,m}}) & \forall p \in P_{ij}\\
    \tau_{oj}^p &= 0 & \forall j \in C, p \in P_{ij}\\
    \sum\limits_{g \in Co} \sum\limits_{p \in P_{ig}} \tau_{ig}^p &= \sum\limits_{j \in Co} \sum\limits_{p \in P_{ji}} \tau_{ji}^p + t_{ji}^{p} + \Delta_{ji}^{p} & \forall i \in C \\
    \tau_{ij}^{p} &\leq T x_{ij}^{p} & \forall i, j \in Co, p \in P_{ij}\\
    \tau_{io}^{p} &\leq T x_{io}^{p} - \Delta_{io}^{p} - t_{io}^{p} x_{io}^{p} & \forall i \in C, p \in P_{ij}\\
    x_{ij}^{p} &\leq y_h  & \forall i, j \in Co, p \in P_{ij}, h \in p\\
    \sum\limits_{h \in S} y_h &\leq \bar{S}\\
    \sum\limits_{j \in C} \sum\limits_{p \in P_{oj}} x_{oj}^{p} &\leq \bar{R}\\
    \tau_{ij}^{p}, q_{ij}^{p} &\geq 0 & \forall p \in P_{ij}\\
    {q^+}_{ij}^{p}, {q^-}_{ij}^{p} &\geq 0 & \forall p \in P_{ij} \\
    s_{ij}^{p}, d_{ij}^{p}, \Delta_{ij}^{p} &\geq 0 & \forall p \in P_{ij} \\
    {\alpha^+}_{ijk}^{p}, {\alpha^-}_{ijk}^{p} &\in \left[0, 1\right] & \forall p \in P_{ij}, k \in B \\
    x_{ij}^{p} &\in \{0,1\} & \forall p \in P_{ij}\\
    y_h &\in \{0,1\} & \forall h \in S \\
    {z^+}_{ijk}^{p}, {z^-}_{ijk}^{p} &\in \{0,1\} & \forall p \in P_{ij}, k \in B
\end{align}

The objective function (90) aims to minimize the total time, similar to the previous models. Constraints (91) ensure that each customer is visited once, and (92) enforce flow conservation. Constraints (93) define the battery SoC as equal to $Q$ when an EV departs from the depot to travel recharge path $p$. Constraints (94) and (95) define the battery SoC when an EV arrives at the first CS on path $p$ after departing from the depot and from a customer node, respectively. \change{The parameter $e^{p,f}_{ij}$ refers to the energy cost of the arc connecting node $i$ to the first CS in path $p$.} Constraints (96) couple the battery level of an EV between two \change{CSs} in a path $p$, if $|p| \geq 2$. The parameter \change{$e^{p,m}_{ij}$} represents the energy arc incident to the CS at the position $m$ in path $p \in P_{ij}$.  Constraint (97) track the EV battery SoC when an EV departs from each customer $i$. Constraints (98) and (99) define that if no EV travels on arc $x_{ij}^p$, then $q_{ij}^p$ and ${q^+}_{ij}^{p}$ is 0. Constraints (100) ensure that if no EV travels on path $x_{ij}^p$, no battery recharge can take on $p$. Constraints (101) ensure that EVs must have sufficient energy to complete the route, and (102) ensure that $q_{ij}^p = 0$ if any CS $h \in p$ is not opened.

The piecewise linear function formulation follows the same principle as the previous models, with the remodeled variables: ${\alpha^+}_{ijk}^\change{{p,m}}$ is the coefficient of each energy segment $k \in B$ associated to battery level ${q^+}_{ij}^\change{{p,m}}$ in the CS $m \in p$, and ${\alpha^-}_{ijk}^\change{{p,m}}$ is the coefficient of energy segment $k \in B$ associated to battery level ${q^-}_{ij}^\change{{p,m}}$; the binary variable $z^+_{ijk}$ is equal to 1 if ${q^+}_{ij}^\change{{p,m}} \geq \sum_{k \in B} \change{E^{'p,m}}_{ij(k-1)}$ and ${z^-}_{ijk}^\change{{p,m}}$ is equal to 1 if ${q^-}_{ij}^\change{{p,m}} \geq \sum_{k \in B} \change{E^{'p,m}}_{ij(k-1)}$ where $\change{E^{'p,m}}_{ijk}$ refers to the $k$-th energy segment value of the piecewise charging function of CS $m \in p$. $\change{T^{'p,m}}_{ijk}$ refers to the $k$-th time segment value of the charging function of CS $m \in p$.  Constraints (103)-(107) and (108)-(112) respectively define the battery charge level and its corresponding charging time when an EV enters and leaves a CS in the $m \in p$ if $|p| \geq 1$. Constraints (113) compute the time spent at recharge path $p$.

Constraints (114) define that the time when an EV departs from the depot is equal to 0, and (115) track the EV time when it departs from each customer $i$. Constraints (116) define that if no EV travels on arc $x_{ij}^p$, then $\tau_{ij}^p$ is 0. Constraints (117) ensure that any route is completed no later than $T$. Constraints (118) ensure that an EV can only travel a recharge path $p$ from $i$ to $j$ if all the \change{CSs} $h \in p$ are opened. The number of opened \change{CSs} is limited to $\bar{S}$ by constraints (119), and (120) limit the number of EV routes. Constraints (121)-(127) define the variable domains.

\subsubsection{Preprocessing of recharge paths}

This Section proposes a \change{preprocessing procedure to identify strictly worse recharge paths} in the third and fourth models that may not appear in an optimal solution and can be excluded from the problem to reduce the search space and improve the performance of the solution process. We will refer here only to recharge paths, as a recharge arc can be treated as a recharge path with only one CS.

\change{We use as a reference to the preprocessing procedure the dominance analysis in recharge paths presented by \cite{froger2019improved} in the context of EVRP models with recharge paths and similar objective function, from which we adapted the following definition: A recharge path $p'$ is said to be dominated by another path $p$ between nodes $i,j ,\in C \cup \{o\}$ if, for every possible SoC at the destination $j$, it is possible for an EV to travel from $i$ to $j$ via path $p$ in a shorter time duration than that of path $p'$ (time-dominance condition), and for every possible time at the destination $j$, it is also possible for an EV to travel from $i$ to $j$ using path $p$ with lower battery consumption compared to path $p'$, reaching $j$ with a higher state of charge (energy-dominance condition). Therefore, if time-dominance and energy-dominance condition is established comparing two paths $p$ and $p'$, than path $p'$ will be strictly worse then path $p$ and may be excluded without compromising the optimal solution}.

\change{For each pair of nodes $i,j \in Co$ in the ELRP, there is a set of potential recharge paths, but not all of them may be part of an optimal solution. Among the potential recharge paths between $(i,j)$, the preprocessing procedure evaluates the time-dominance and energy-dominance condition to find the paths that are strictly worse with respect to energy and time consumption. A feasible path that connects directly $i$ to $j$ without CS cannot be dominated by other paths.}

\change{Algorithm 1 illustrates the procedure. It starts by analyzing the paths with $k = 1$ CS. For each feasible path $p$ between each pair of nodes $(i,j)$ we calculate the total travel time cost $ct = t^p_{ij}$, the total travel energy cost $ce = e^p_{ij}$, the first and last arc energy cost $cf = e_{if}$ and $cl = e_{lj}$ where node $f$ is the first CS and node $l$ is the last CS in the path. In addition to influencing the total energy cost, $cf$ also influences the energy feasibility of the path, as a lower $cf$ means a lower SoC will be necessary to allow an EV to reach the first CS in $p$. The last arc $cl$ is critic because it influences the SoC that the EV can achieve at the end of the path. The recharge speed related to the CS type is represented by $sp$. Subsequently, the feasible paths are ordered in $p_{ij}$ from lowest to highest $ct$ and a pairwise comparison ($p$, $p'$) is performed in sequence, starting from the first to the last path in $p_{ij}$. The comparison evaluates the time-dominance and energy-dominance conditions of the paths according to the relationship between the station speed of each path. If $ct < ct' $, $cf < cf' $, $cl < cl' $ and $sp \geq sp'$ then the time-dominance and energy-dominance condition is identified and $p'$ is strictly worse than $p$. We call this condition Case 1.}

\change{If the recharge speed of station $s$ in $p$ is slower than the station $s'$ in $p'$ ($sp < sp'$), we call this Case 2, and we should consider in the comparison the maximum difference in charging time and charged energy between CSs in each path:
$\bar{\Delta t} = max(\Delta t_s - \Delta t_{s'})$ and $\bar{\Delta q} = max(\Delta q_s - \Delta q_{s'})$. These upper bounds help determine the conditions under which time and energy dominance between paths can be asserted in Case 2. In this case, if
$ct < ct' $, $ce < ce'$, $cf < cf'$, $cl < cl'$, $ct + \bar{\Delta t} < ct'$ and $ce + \bar{\Delta q} < ce'$, then the time-dominance and energy-dominance condition is identified and $p'$ is strictly worse than $p$. In Appendix A we present a detailed description of the time-dominance and energy-dominance evaluation of Case 1 and Case 2. After the pairwise comparison among all paths, the paths identified in the time and energy dominance condition are excluded, and the procedure is repeated from line 2, increasing $k$ by $1$ until $k \leq \bar{S}$ or until there are no more dominated paths to be excluded. It is worth noting that due to the triangular inequality, the total energy cost and travel time cost will always increase when one station is added to a recharge path. The algorithm returns the set $Pd$ that contains all non-dominated recharge paths.}

\begin{algorithm}[htbp]
\caption{Identify and exclude dominated recharge paths}\label{alg:cap}
\begin{algorithmic}[1]
\renewcommand{\algorithmicrequire}{\textbf{Input:}} 
\renewcommand{\algorithmicensure}{\textbf{Output:}} 
\Require Sets $C_o , S, Pd=\{\} $
\Ensure Set $Pd$ with all non-dominated $p$ between each nodes pair $(i,j)$
\State $k \gets 1$
\State $p_{ij} \gets$ all feasible $p$ with k stations between nodes $(i,j)$
\State $Pd \gets$ all $p_{ij}$
\State \textbf{for} $p$ in $p_{ij}$
\State \qquad Perform a pairwise comparison between each $p$
\State \qquad Classify comparison in Case 1 or Case 2
\State \qquad Evaluate time-dominance and energy-dominance
\If{ time-dominance and energy-dominance is established}
\State Exclude $p'$ from $Pd$
\EndIf
\State update $k \gets  k + 1$
\If{$k \leq |\bar{S}|$ and there are still dominated paths to be excluded}
\State \textbf{go to} Line 2
\EndIf
\State \textbf{return} Set $Pd$ containing all non-dominated recharge paths
\end{algorithmic}
\end{algorithm}

\section{Computational Experiments}

In this section, the results of computational experiments are presented and discussed. The models were implemented using \change{Julia Programming Language (version 1.10.7), JuMP Modeling Language (version 1.23.6), and Gurobi Solver (version 12.0). The results of the four formulations proposed for ELRP-NLMS were performed under the same conditions using a computer with two Intel (R) Xeon (R) E5-2640 v4 2.40GHz CPUs with 20 physical cores and 40 threads, with 256GB of RAM, using up to one thread, and with a runtime limit of 10,800 seconds (3 hours) for each instance.}

\subsection{ELRP-NLMS Instance Generation}

Instances presented by \cite{montoya2017electric} for the Electric Vehicle Routing Problem with nonlinear charging function (EVRP-NL) were extended to the ELRP. These instances were generated based on a geographic space of $120 \times 120$ km, representing a semi-urban operation using real data for EVs and charging function settings. The EVs have a consumption rate of 0.125 kWh/km and a battery capacity of 16 kWh. The maximum duration of time for each route is 10 hours. \change{The energy consumption on an arc $e_{ij}$ is calculated based on the EV’s consumption rate multiplied by the arc’s length, following a classical approach in the literature. \cite{montoya2017electric} considered two levels of availability of the charging infrastructure in EVRP-NL instances: low and high. To ensure feasibility, the number of CSs was selected as a proportion of the number of customers for each combination of customer count and infrastructure availability. The CSs are located either randomly or through a simple p-median heuristic, which starts with randomly generated CS locations and iteratively adjusts them to minimize the total distance to customers. Three types of CS are included: slow, moderate, and fast, with their types randomly selected using a uniform probability distribution. More details on EVRP-NL instances are available at \cite{montoya2017electric}.}

\change{The extension of EVRP model instances to ELRP was executed using the following methodology. In addition to the original $k$ CSs, other $k$ customer nodes are randomly selected to serve as CS candidates.  In this way, the ELRP instance now has $2k$ candidates for CSs, keeping the same limit of $k$ CSs from the original instance that can be opened. The number of CSs is doubled to ensure that for each CS opening decision, the model has at least two candidates to choose from.} Regarding the type of CSs, new candidates for CSs will always be set as the slowest CS type from the original instance. This means that if the original instance has two CS types, one fast and one slow, both new potential stations will be of the slow type. This approach is designed to prevent bias when comparing the results of the ELRP with the EVRP. This methodology aligns with the approach taken by \cite{schiffer2017electric}, who also extended EVRP instances to ELRP by transforming customer nodes into potential CS candidates. However, a key difference is that Schiffer et al. consider all customer nodes as CS candidates. It is important to note that the decision to open a CS is not always feasible at every location, as it depends on factors such as parking availability, energy grid conditions, and power capacity. \change{Therefore, by limiting the number of customer nodes to CS candidates and doubling $k$, our approach provides a reasonable number of CS candidates while maintaining a more realistic framework.}

A set of 80 instances were adapted to the ELRP-NLMS model. The set of instances can be grouped into four subsets according to the size of the instance: 20 instances with ten customers and 4-6 \change{CSs}; 20 instances with 20 customers and 6-8 \change{CSs}; 20 instances with 40 customers and 10-16 \change{CSs}; and 20 instances with 80 customers and 16-24 \change{CSs}. As the first and second model formulations use $\beta$ parameter to replicate the charge stations, we use a similar method as presented in \cite{montoya2017electric} to define this parameter: Starting with $\beta = 0$, we try to solve the instance with a time limit of 3 hours and then, at each subsequent iteration, we solve the model with $\beta = \beta + 1$. The procedure stops when the time limit is reached and no improvement in the optimal or best solution is found.

\subsection{\change{Result Analysis}}

Table 2 provides a summary of the performance of the four formulations grouped by instance sizes, in which M1, M2, M3, and M4 refer to the first, second, third and fourth model formulations, respectively. To evaluate each group of instances for each model formulation, we computed several metrics: the average MIP Gap (Avg Gap), the number of optimally solved instances (\#Opt), the number of instances with feasible solutions that were not solved to optimality (\#Feas), \change{the average runtime for solved instances (Avg Time), the number of instances that achieved the best objective function value and the best MIP Gap among the four formulations (\#BestObj and \#BestGap), the average Linear Relaxation (Avg LR) and the minimum Gap (MinGap). To measure the MIP Gap, we use the optimality Gap as reported by the Gurobi solver, computed by $(|ObjBound - ObjVal| / ObjVal) \times 100$ where $ObjBound$ and $ObjVal$ refer respectively to the best objective relaxed bound and objective value found by the solver. The \#BestObj and \#BestGap considered instances with at least two feasible solutions to make the comparison. The Avg LR was calculated by removing the integrality constraints of each formulation variable and solving each instance to optimality. The detailed results for the 80 instances and each formulation are presented in Appendix B.}

\change{Preliminary tests were conducted to asses the impact of the preprocessing strategies outlined in Section 2.2. Valid inequalities (38) and (39) have reduced the gap in M1 and M2 by an average of 1.2\% in experiments with test instances. The Subtour Elimination Constraints have reduced the average gap from 37.8\% to 28.9\% in tests with 54 instances in M1 but did not show improvements in M2. Therefore, in the final experiments we used the SEC only in M1. The application of the Subtour Elimination Constraints was implemented in two stages: initially, the integrality constraints of the models were removed, and they were solved with a cutting plane algorithm. In the second stage, the integrality of the model variables was restored, and the SECs were added using the Gurobi solver callback function, resulting in a branch-and-cut. Table 3 presents the average gap and the average time for the SEC evaluation with 54 instances in M1 and M2.}

\begin{table}[htbp]
  \centering
    \caption{General results per instance group.}
    \small
    \addtolength{\tabcolsep}{-1pt}
    \begin{tabular}{cccccrccccc}
    \multicolumn{5}{p{21.855em}}{10-customers instances} &       & \multicolumn{5}{p{20.715em}}{20-customers instances} \\
\cmidrule{1-5}\cmidrule{7-11}          & M1    & M2    & M3    & M4    &       &       & M1    & M2    & M3    & M4 \\
\cmidrule{1-5}\cmidrule{7-11}    \textbf{Avg Gap} & 12.2\% & 0.0\% & 0.0\% & 0.0\% &       & \textbf{Avg Gap} & 30.4\% & 13.3\% & 12.4\% & 14.5\% \\
    \textbf{\#Opt} & 11    & 20    & 20    & 20    &       & \textbf{\#Opt} & 0     & 8     & 7     & 6 \\
    \textbf{\#Feas} & 9     & 0     & 0     & 0     &       & \textbf{\#Feas} & 20    & 12    & 13    & 14 \\
    \textbf{Avg Time} & 5703  & 444   & 348   & 779   &       & \textbf{Avg Time} & 10800 & 7104  & 7827  & 9231 \\
    \textbf{\#BestObj} & 19    & 20    & 20    & 20    &       & \textbf{\#BestObj} & 10    & 15    & 16    & 15 \\
    \textbf{\#BestGap} & 11    & 20    & 20    & 20    &       & \textbf{\#BestGap} & 0     & 11    & 14    & 6 \\
    \textbf{Avg LR} & 3.544 & 5.654 & 6.087 & 6.089 &       & \textbf{Avg LR} & 5.520 & 9.133 & 9.802 & 9.799 \\
    \textbf{MinGap} & 0.0\% & 0.0\% & 0.0\% & 0.0\% &       & \textbf{MinGap} & 14.1\% & 0.0\% & 0.0\% & 0.0\% \\
\cmidrule{1-5}\cmidrule{7-11}          &       &       &       &       &       &       &       &       &       &  \\
    \multicolumn{5}{p{21.855em}}{40-customers instances} &       & \multicolumn{5}{p{20.715em}}{80-customers instances} \\
\cmidrule{1-5}\cmidrule{7-11}          & M1    & M2    & M3    & M4    &       &       & M1    & M2    & M3    & M4 \\
\cmidrule{1-5}\cmidrule{7-11}    \textbf{Avg Gap} & \multicolumn{1}{r}{46.3\%} & \multicolumn{1}{r}{22.6\%} & \multicolumn{1}{r}{25.2\%} & \multicolumn{1}{r}{37.8\%} &       & \textbf{Avg Gap} & \multicolumn{1}{r}{48.8\%} & \multicolumn{1}{r}{29.2\%} & \multicolumn{1}{r}{44.0\%} & \multicolumn{1}{r}{100.0\%} \\
    \textbf{\#Opt} & 0     & 0     & 0     & 0     &       & \textbf{\#Opt} & 0     & 0     & 0     & 0 \\
    \textbf{\#Feas} & 16    & 13    & 13    & 13    &       & \textbf{\#Feas} & 17    & 7     & 13    & 4 \\
    \textbf{Avg Time} & 10800 & 10800 & 10800 & 10800 &       & \textbf{Avg Time} & 10800 & 10800 & 10800 & 10800 \\
    \textbf{\#BestObj} & 4     & 9     & 2     & 0     &       & \textbf{\#BestObj} & 5     & 5     & 2     & 0 \\
    \textbf{\#BestGap} & 0     & 9     & 4     & 1     &       & \textbf{\#BestGap} & 1     & 7     & 4     & 0 \\
    \textbf{Avg LR} & 7.546 & 14.255 & 15.390 & 15.395 &       & \textbf{Avg LR} & 11.683 & 21.049 & 21.920 & 21.920 \\
    \textbf{MinGap} & 39.2\% & 11.0\% & 11.9\% & 20.7\% &       & \textbf{MinGap} & 43.3\% & 20.7\% & 36.2\% & 100.0\% \\
\cmidrule{1-5}\cmidrule{7-11}    \end{tabular}%
  \label{tab:addlabel}%
\end{table}%

\begin{table}[htbp]
  \centering
  \caption{SEC evaluation with 54 instances in M1 and M2}
  \small
    \addtolength{\tabcolsep}{-1pt}
    \begin{tabular}{ccrc}
    \toprule
    \textbf{Model} & \textbf{Subtour Elimination Constraints} & \multicolumn{1}{c}{\textbf{Avg Gap}} & \textbf{Avg Time} \\
    \midrule
    NodeModel & SEC   & 28.88\% & 8758 \\
    NodeModel & None  & 37.85\% & 8956 \\
    ArcModel & SEC   & 12.47\% & 5529 \\
    ArcModel & None  & 12.46\% & 5497 \\
    \bottomrule
    \end{tabular}%
  \label{tab:addlabel}%
\end{table}%

\change{Considering the results of the four models, feasible solutions were obtained for 74 out of the 80 instances. Figure 2 illustrates the overall performance per model type, showing the total number of instances optimally solved and the number of instances that achieved the \#BestObj and \#BestGap  among the four models. The results highlight that formulation M2 shows the best performance among the four formulations, achieving 28 optimal solutions, 49 \#BestObj and 47 \#BestGap out of the 80 instances. Following M2, formulation M3 stands out, with 27 optimal solutions, 40 \#BestObj and 42 \#BestGap. Formulation M4 obtained 26 optimal solutions and 27 \#BestGap, outperforming M1 in these metrics, which attained 11 optimal solutions and 12 \#BestGap. However, M1 achieved 38 \#BestObj, surpassing M4 which obtained 35 \#BestObj.}

\change{It is worth noting that M4 performed poorly in the group of instances with 80 customers, obtaining feasible solutions for only 4 instances with an average Gap of 100\%, and that M1 obtained a greater number of feasible solutions in the groups of instances with 40 and 80 customers. It is noted that formulation M4 has the largest number of variables and constraints and that the limitation of the computer's processor to up to 1 thread in the experiments can make it harder the solution of larger instances in these formulations. In future research we propose new experiments with more powerful computers in larger instances.}

\begin{figure}[htbp]
  \centerline{\includegraphics[width=12cm]{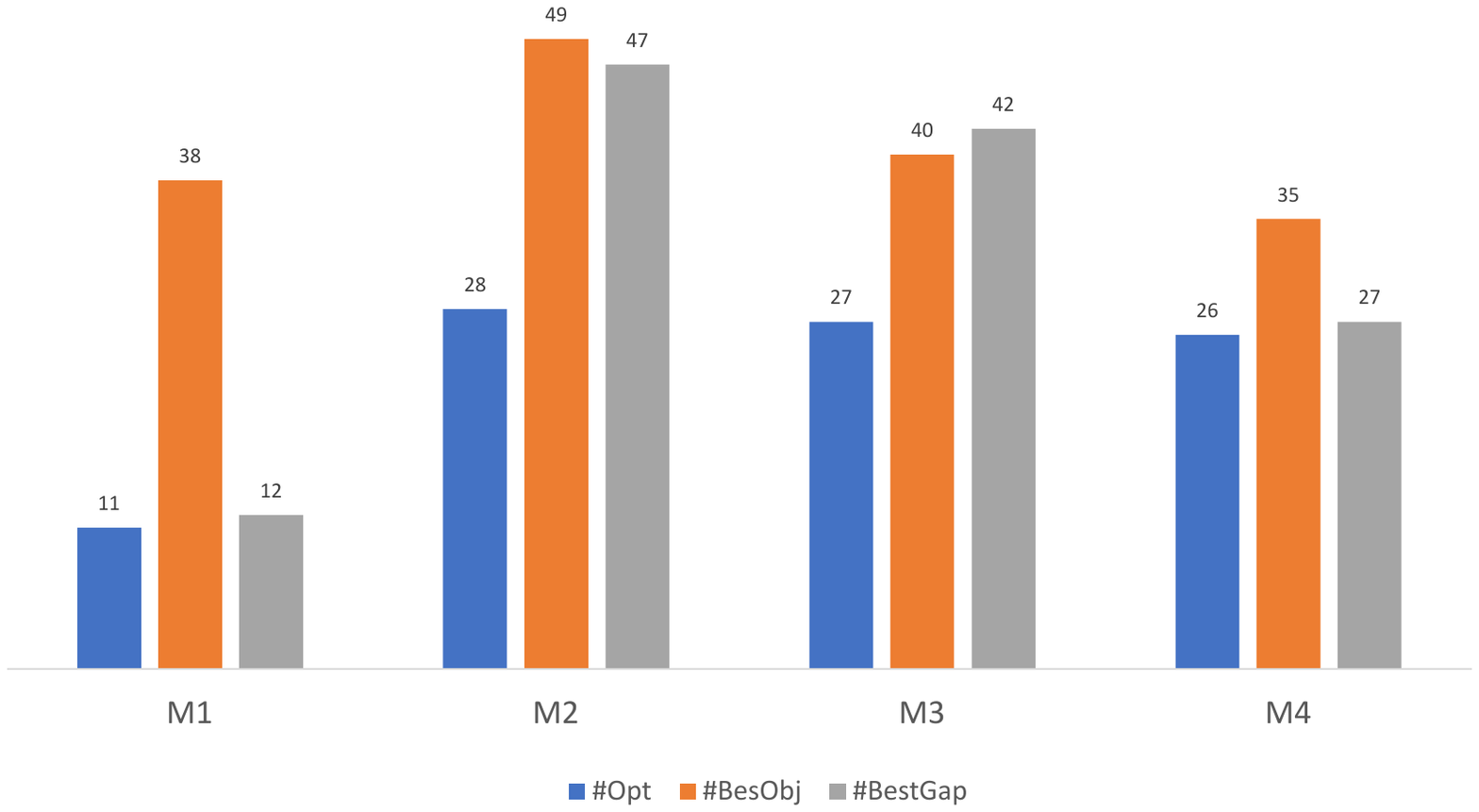}}
\caption{General results per model type.\label{fig:3}}
\end{figure}

\change{The Linear Relaxation (LR) is a useful metric as it helps to find bounds on the optimal solution of MIP. We ran the four models in a relaxed version to evaluate the quality of bounds provided by each formulation. The results for the LR of each formulation for all 80 instances are presented in Appendix C. The results indicate that all three new proposed formulations achieve significantly better bounds than the first formulation, with emphasis on M3 and M4, followed by M2. Table 4 presents the integrality Gap of the four formulations in the 20 customer instances, comparing the values of the objective function in the optimal solution of the MIP and the relaxed version. The results show that the average integrality Gap of M1 is 70.3\%, significantly higher than M2 with 52.7\% and M3 and M4 with 48.8\%.}

\begin{table}[htbp]
  \centering
  \caption{Integrality Gap of 10 customer instances}
  \scriptsize
    \begin{tabular}{rrrcccccccc}
    \toprule
          &       &       & \multicolumn{2}{c}{\textbf{M1}} & \multicolumn{2}{c}{\textbf{M2}} & \multicolumn{2}{c}{\textbf{M3}} & \multicolumn{2}{c}{\textbf{M4}} \\
\cmidrule{4-11}    \multicolumn{1}{c}{\textbf{\#}} & \multicolumn{1}{c}{\textbf{Instance}} & \textbf{MIP \newline{}Opt} & \multicolumn{1}{p{3.57em}}{\textbf{RL}} & \multicolumn{1}{p{3.57em}}{\textbf{Int. Gap}} & \multicolumn{1}{p{3.57em}}{\textbf{RL}} & \multicolumn{1}{p{3.57em}}{\textbf{Int. Gap}} & \multicolumn{1}{p{3.57em}}{\textbf{RL}} & \multicolumn{1}{p{3.57em}}{\textbf{Int. Gap}} & \multicolumn{1}{p{3.57em}}{\textbf{RL}} & \multicolumn{1}{p{3.57em}}{\textbf{Int. Gap}} \\
    \midrule
    \multicolumn{1}{c}{1} & \multicolumn{1}{c}{tc0c10s2cf1-p4} & 19.753 & 5.808 & 70.6\% & 7.506 & 62.0\% & 7.291 & 63.1\% & 7.348 & 62.8\% \\
    \multicolumn{1}{c}{2} & \multicolumn{1}{c}{tc0c10s2ct1-p4} & 11.378 & 5.057 & 55.6\% & 6.894 & 39.4\% & 7.322 & 35.6\% & 7.310 & 35.8\% \\
    \multicolumn{1}{c}{3} & \multicolumn{1}{c}{tc0c10s3cf1-p6} & 10.921 & 5.223 & 52.2\% & 6.467 & 40.8\% & 7.205 & 34.0\% & 7.205 & 34.0\% \\
    \multicolumn{1}{c}{4} & \multicolumn{1}{c}{tc0c10s3ct1-p6} & 10.536 & 4.381 & 58.4\% & 6.235 & 40.8\% & 7.218 & 31.5\% & 7.218 & 31.5\% \\
    \multicolumn{1}{c}{5} & \multicolumn{1}{c}{tc1c10s2cf2-p4} & 9.034 & 3.864 & 57.2\% & 4.765 & 47.3\% & 5.571 & 38.3\% & 5.571 & 38.3\% \\
    \multicolumn{1}{c}{6} & \multicolumn{1}{c}{tc1c10s2cf3-p4} & 12.583 & 2.939 & 76.6\% & 6.982 & 44.5\% & 6.628 & 47.3\% & 6.664 & 47.0\% \\
    \multicolumn{1}{c}{7} & \multicolumn{1}{c}{tc1c10s2cf4-p4} & 16.097 & 5.253 & 67.4\% & 6.928 & 57.0\% & 7.226 & 55.1\% & 7.226 & 55.1\% \\
    \multicolumn{1}{c}{8} & \multicolumn{1}{c}{tc1c10s2ct2-p4} & 9.199 & 3.435 & 62.7\% & 4.652 & 49.4\% & 5.571 & 39.4\% & 5.571 & 39.4\% \\
    \multicolumn{1}{c}{9} & \multicolumn{1}{c}{tc1c10s2ct3-p4} & 12.307 & 2.743 & 77.7\% & 5.896 & 52.1\% & 6.481 & 47.3\% & 6.481 & 47.3\% \\
    \multicolumn{1}{c}{10} & \multicolumn{1}{c}{tc1c10s2ct4-p4} & 13.826 & 5.192 & 62.4\% & 6.722 & 51.4\% & 7.224 & 47.8\% & 7.224 & 47.8\% \\
    \multicolumn{1}{c}{11} & \multicolumn{1}{c}{tc1c10s3cf2-p6} & 9.034 & 4.230 & 53.2\% & 5.293 & 41.4\% & 5.571 & 38.3\% & 5.571 & 38.3\% \\
    \multicolumn{1}{c}{12} & \multicolumn{1}{c}{tc1c10s3cf3-p6} & 12.583 & 2.452 & 80.5\% & 5.741 & 54.4\% & 6.480 & 48.5\% & 6.480 & 48.5\% \\
    \multicolumn{1}{c}{13} & \multicolumn{1}{c}{tc1c10s3cf4-p6} & 14.902 & 4.917 & 67.0\% & 6.727 & 54.9\% & 7.348 & 50.7\% & 7.342 & 50.7\% \\
    \multicolumn{1}{c}{14} & \multicolumn{1}{c}{tc1c10s3ct2-p6} & 9.199 & 3.992 & 56.6\% & 5.022 & 45.4\% & 5.571 & 39.4\% & 5.571 & 39.4\% \\
    \multicolumn{1}{c}{15} & \multicolumn{1}{c}{tc1c10s3ct3-p6} & 12.932 & 2.301 & 82.2\% & 5.727 & 55.7\% & 6.480 & 49.9\% & 6.480 & 49.9\% \\
    \multicolumn{1}{c}{16} & \multicolumn{1}{c}{tc1c10s3ct4-p6} & 13.205 & 5.290 & 59.9\% & 6.688 & 49.4\% & 7.355 & 44.3\% & 7.319 & 44.6\% \\
    \multicolumn{1}{c}{17} & \multicolumn{1}{c}{tc2c10s2cf0-p4} & 11.345 & 1.034 & 90.9\% & 3.728 & 67.1\% & 3.797 & 66.5\% & 3.797 & 66.5\% \\
    \multicolumn{1}{c}{18} & \multicolumn{1}{c}{tc2c10s2ct0-p4} & 11.345 & 0.886 & 92.2\% & 3.717 & 67.2\% & 3.797 & 66.5\% & 3.797 & 66.5\% \\
    \multicolumn{1}{c}{19} & \multicolumn{1}{c}{tc2c10s3cf0-p6} & 11.285 & 1.020 & 91.0\% & 3.709 & 67.1\% & 3.797 & 66.4\% & 3.797 & 66.4\% \\
    \multicolumn{1}{c}{20} & \multicolumn{1}{c}{tc2c10s3ct0-p6} & 11.285 & 0.873 & 92.3\% & 3.689 & 67.3\% & 3.797 & 66.4\% & 3.797 & 66.4\% \\
    \midrule
          &       & Average &       & 70.3\% &       & 52.7\% &       & 48.8\% &       & 48.8\% \\\bottomrule
    \end{tabular}%
  \label{tab:addlabel}%
\end{table}%

\change{Figure 2 presents a Box Plot chart with the Gap distribution for each formulation and each instance group, emphasizing the reduced Gap achieved by the three newly proposed formulations. For the Gap distribution, we considered only the instances with feasible solutions in the four formulations to make the comparison balanced, with the exception of the 80-customer instances, in which we disregarded M4, which provided solutions for only four instances with a Gap of 100\%. In the instances with ten customers, M2, M3 and M4 achieved the optimal solution in all 20 instances, while M1 achieved the optimal only in nine instances, resulting in an average Gap of 12.2\%.  For the group of 20-customer instances, M3 found the lowest average Gap, 13.3\%,  but M2 optimally solved eight instances, one more than M3. In the 40-customer instances, M2 achieved the best average and minimum Gap, 22.6\% and 11.0\% respectively. In the 80-customer instances,  M2 also achieved the best average and minimum Gap, 29.2\% and 20.7\% respectively. When comparing all solutions obtained in M1 and M2, the average Gap was reduced from 29.1\% to 11.9\%.
In Appendix D, we consolidate the best results from each instance in the experiments, highlighting the best objective function and the best lower bound obtained by the Gurobi solver in the four formulations, resulting in the achievement of even smaller Gaps. In this case, the minimum Gap in the 40-customer instances was reduced to 10.2\%, while in the 80-customer group it was reduced to 17.8\%.} 

\begin{figure}[htbp]
  \centerline{\includegraphics[width=15cm]{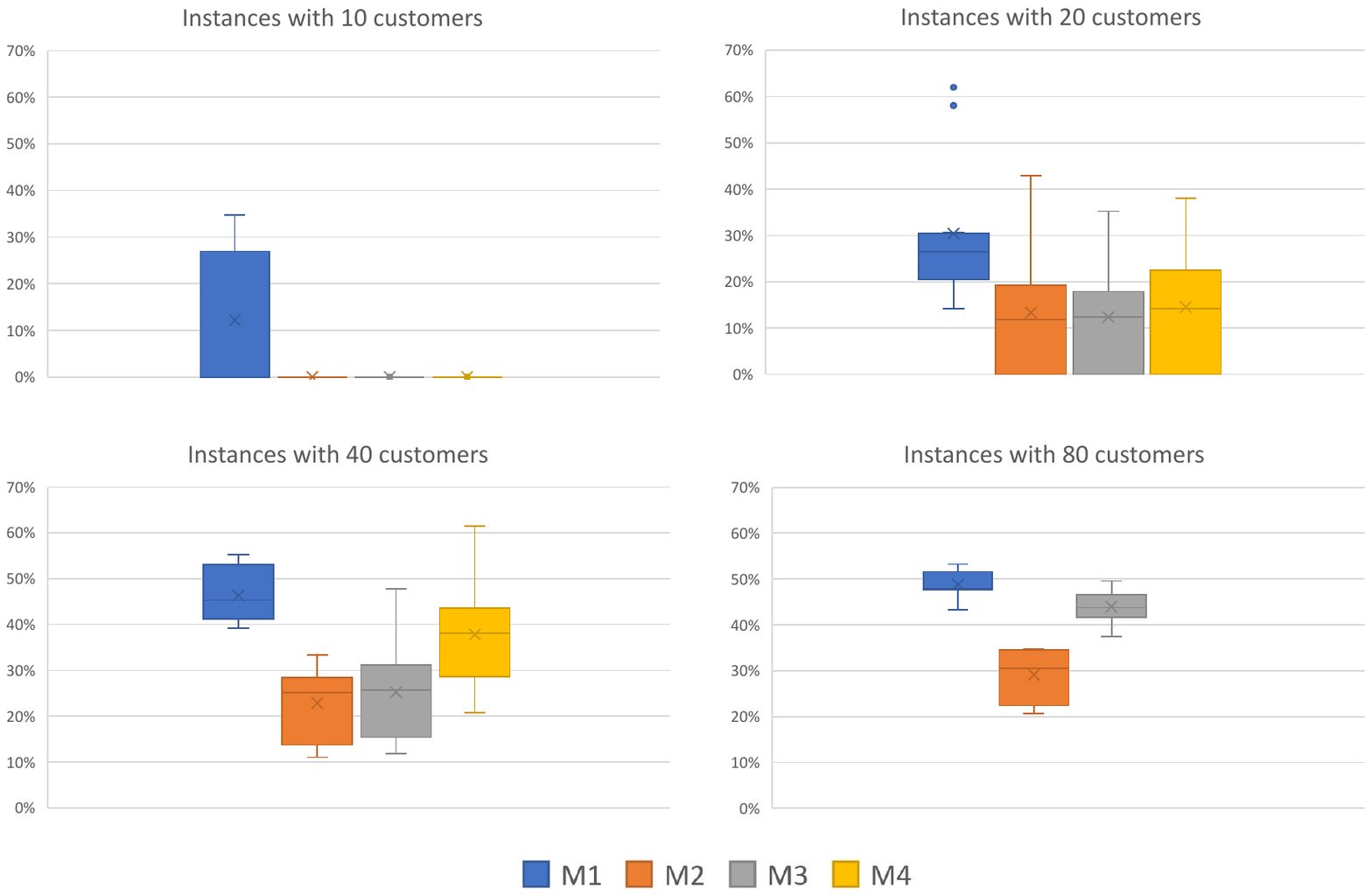}}
\caption{ Optimality Gap comparison per instance group.\label{fig:2}}
\end{figure}

\subsection{Effects of nonlinear charging function on ELRP decisions}

\change{The ELRP simultaneously integrates EV routing and CS location decisions to optimize the performance of a logistic network. In this model, each customer must be serviced exactly once by a single EV, and the number of CSs that can be installed is limited to an upper bound $\bar{S}$. Although the actual charging process of EV batteries has a nonlinear behavior (Pelletier et al., 2017), most ELRP models in the literature simplify this by assuming linear charging. Therefore, this section investigates the implications of adopting linear and nonlinear charging assumptions in ELRP decision-making}.

\change{To evaluate the impact of the nonlinear charging function on EV location-routing decisions, we defined the  electric location-routing model with linear charging and multiple charging station types (ELRP-L-MS). This model is analogous to the ELRP-NLMS defined in Section 2, with the key difference being the replacement of the nonlinear charging function with a linear one. Figure 4 illustrates the linear charging function, represented by the red dotted line, in comparison with the nonlinear piecewise approximation function. Furthermore, we executed the ELRP model while maintaining the route and CS opening decisions from the ELRP-L-MS, but applied the nonlinear charging function (referred to as ELRP-NL-LD). This approach allows us to assess the impact of the nonlinear charging function on the objective function results.}

\begin{figure}[htbp]
  \centerline{\includegraphics[width=10cm]{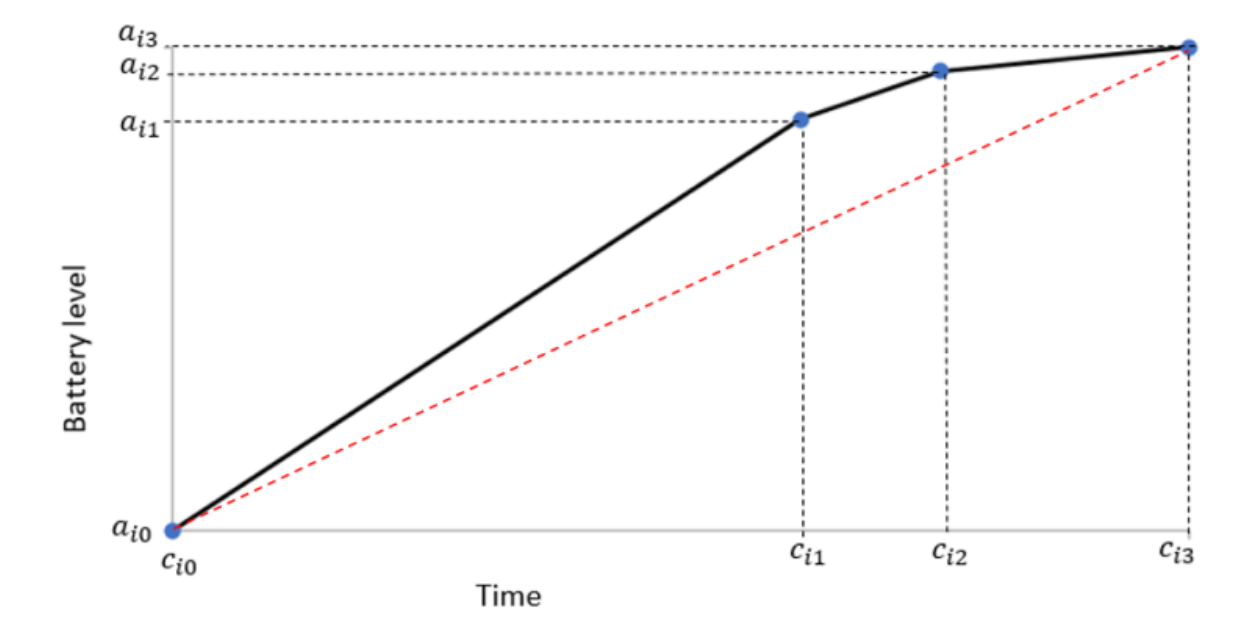}}
\caption{Linear charge function and non-linear piecewise approximation charging function.\\ (Adapted from Montoya et al., 2017)\label{fig:4}}
\end{figure}

Table 5 presents the ELRP-NLMS and the ELRP-L-MS detailed results for the twenty instances with ten customers optimally solved. The first two columns identify the instances; the following six columns present the results of the ELRP-NLMS and ELRP-L-MS models with the \change{objective function}, the number of EV routes, and the CS opening decisions. The last two columns present the \change{objective function} of the ELRP-NL-LD and its Gap compared with the ELRP-NLMS (NL-G).\\

\begin{table}[htbp]
  \centering
  \caption{Results for the 10-customer instances by models ELRP-NLMS, ELRP-L-MS and ELRP-NL-LD}
  \scriptsize
    \begin{tabular}{cccrcccrcccrr}
    \toprule
        &     &     & \multicolumn{3}{c}{\textbf{        ELRP-NLMS}} &     & \multicolumn{3}{c}{\textbf{  ELRP-L-MS}} &     & \multicolumn{2}{c}{\textbf{ELRP-NL-LD}} \\
\cmidrule{4-6}\cmidrule{8-10}\cmidrule{12-13}    \#  & \textbf{Instance} &     & \multicolumn{1}{c}{\textbf{Obj}} & \multicolumn{1}{c}{\textbf{\#EV \ routes}} & \multicolumn{1}{c}{\textbf{Open\newline{} stations}} &     & \multicolumn{1}{c}{\textbf{Obj}} & \multicolumn{1}{c}{\textbf{\#EV routes}} & \multicolumn{1}{c}{\textbf{Open\newline{} stations}} &     & \multicolumn{1}{c}{\textbf{Obj*}} & \multicolumn{1}{p{4.25em}}{\textbf{NL-G}} \\
\cmidrule{1-2}\cmidrule{4-6}\cmidrule{8-10}\cmidrule{12-13}    1   & tc0c10s2cf1-p4 &     & 19.753 & 3   & \textbf{[12, 13]} &     & 20.384 & 3   & \textbf{[8, 13]} &     & \textbf{20.111} & \textbf{1.8\%} \\
    2   & tc0c10s2ct1-p4 &     & 11.378 & 2   & [8, 12] &     & 11.563 & 2   & [8, 12] &     & \textbf{11.497} & \textbf{1.0\%} \\
    3   & tc0c10s3cf1-p6 &     & 10.921 & 2   & \textbf{[2, 3, 12]} &     & 11.105 & 2   & \textbf{[2, 6, 12]} &     & 10.921 & 0.0\% \\
    4   & tc0c10s3ct1-p6 &     & 10.536 & 2   & [2, 3, 14] &     & 10.705 & 2   & [2, 3, 14] &     & 10.536 & 0.0\% \\
    5   & tc1c10s2cf2-p4 &     & 9.034 & 3   & [12] &     & 9.140 & 3   & [12] &     & \textbf{9.122} & \textbf{1.0\%} \\
    6   & tc1c10s2cf3-p4 &     & 12.583 & \textbf{2} & [5, 13] &     & 16.060 & \textbf{3} & [5, 13] &     & \textbf{15.953} & \textbf{26.8\%} \\
    7   & tc1c10s2cf4-p4 &     & 16.097 & 3   & [12, 13] &     & 16.430 & 3   & [12, 13] &     & 16.097 & 0.0\% \\
    8   & tc1c10s2ct2-p4 &     & 9.199 & 3   & \textbf{[4, 12]} &     & 10.621 & 3   & \textbf{[9, 12]} &     & \textbf{10.477} & \textbf{13.9\%} \\
    9   & tc1c10s2ct3-p4 &     & 12.307 & 2   & \textbf{[5, 13]} &     & 13.607 & 2   & \textbf{[12, 13]} &     & \textbf{13.171} & \textbf{7.0\%} \\
    10  & tc1c10s2ct4-p4 &     & 13.826 & 2   & [12, 13] &     & 14.174 & 2   & [12, 13] &     & 13.826 & 0.0\% \\
    11  & tc1c10s3cf2-p6 &     & 9.034 & 3   & [12] &     & 9.140 & 3   & [12] &     & 9.034 & 0.0\% \\
    12  & tc1c10s3cf3-p6 &     & 12.583 & 2   & \textbf{[7, 11, 13]} &     & 13.827 & 2   & \textbf{[7, 12, 13]} &     & \textbf{13.267} & \textbf{5.4\%} \\
    13  & tc1c10s3cf4-p6 &     & 14.902 & 3   & [12, 13, 14] &     & 15.178 & 3   & [12, 13, 14] &     & 15.022 & \textbf{0.8\%} \\
    14  & tc1c10s3ct2-p6 &     & 9.199 & 3   & \textbf{[4, 5, 12]} &     & 10.804 & 3   & \textbf{[4, 12, 13]} &     & \textbf{10.572} & \textbf{14.9\%} \\
    15  & tc1c10s3ct3-p6 &     & 12.932 & 2   & \textbf{[6, 7, 13]} &     & 13.495 & 2   & \textbf{[6, 11, 13]} &     & \textbf{12.943} & \textbf{0.1\%} \\
    16  & tc1c10s3ct4-p6 &     & 13.205 & 2   & \textbf{[12, 13, 14]} &     & 13.713 & 2   & \textbf{[8, 13, 14]} &     & \textbf{13.321} & \textbf{0.9\%} \\
    17  & tc2c10s2cf0-p4 &     & 11.345 & 3   & [2, 6] &     & 11.770 & 3   & [2, 6] &     & \textbf{11.404} & \textbf{0.5\%} \\
    18  & tc2c10s2ct0-p4 &     & 11.345 & 3   & [2, 6] &     & 11.770 & 3   & [2, 6] &     & \textbf{11.404} & \textbf{0.5\%} \\
    19  & tc2c10s3cf0-p6 &     & 11.285 & 3   & [2, 4, 9] &     & 11.638 & 3   & [2, 4, 9] &     & 11.285 & 0.0\% \\
    20  & tc2c10s3ct0-p6 &     & 11.285 & 3   & [2, 4, 9] &     & 11.638 & 3   & [2, 4, 9] &     & \textbf{11.366} & \textbf{0.7\%} \\
    \bottomrule
    \end{tabular}%
  \label{tab:addlabel}%
\end{table}%

The results in Table 5 demonstrate that including the nonlinear charging function in ELRP models significantly impacts strategic and operational decisions and the objective function \change{values} compared to ELRP models that assume linear charging. Regarding the CS sitting decisions, there were differences in the CS location solutions between ELRP-NLMS and ELRP-L-MS in 8 of the 20 instances. These solutions are printed in bold letters in Table 3. Facility location is a strategic decision that usually involves high costs and increased complexity to further changes after the decision has been made. Comparing the number of EV routes, we can observe that in Instance 6, the ELRP-NLMS model needed two EV routes to attend to all customers, one less than the ELRP-L-MS, impacting the EV fleet size decisions. Figure 5 presents the detailed routing and CS location \change{decisions} for Instance 12 in ELRP-L-MS and ELRP-NLMS. As can be seen, in the ELRP-NLMS model, the \change{CSs} at nodes 7, 11, and 12 were opened, while in the ELRP-L-MS, the \change{CSs} 7, 12, and 13 were opened, also changing the EV routes and impacting the objective function \change{value}.

Regarding the objective function \change{results}, the time cost has increased in the ELRP-NL-LD in 14 of the 20 instances, with an average NL-G increase of 3.8\%. The NL-G elucidates the actual \change{objective function} Gap a network manager should have when assuming a linear charging function in his planning model. The higher NL-G is 26.8\% in Instance 6, followed by instances 14 and 8 with 14.9\% and 13.8\%, respectively. In seven instances that kept the same CS opening decisions, the \change{objective function results} has increased due to routing decision changes of the ELRP-L-MS. This result demonstrates that even if the CS sitting decisions are not changed, assuming linear charging could also decrease the \change{objective function} performance of the ELRP model by finding worse routing decisions to the network.

\begin{figure}[htbp]
  \centerline{\includegraphics[width=15cm]{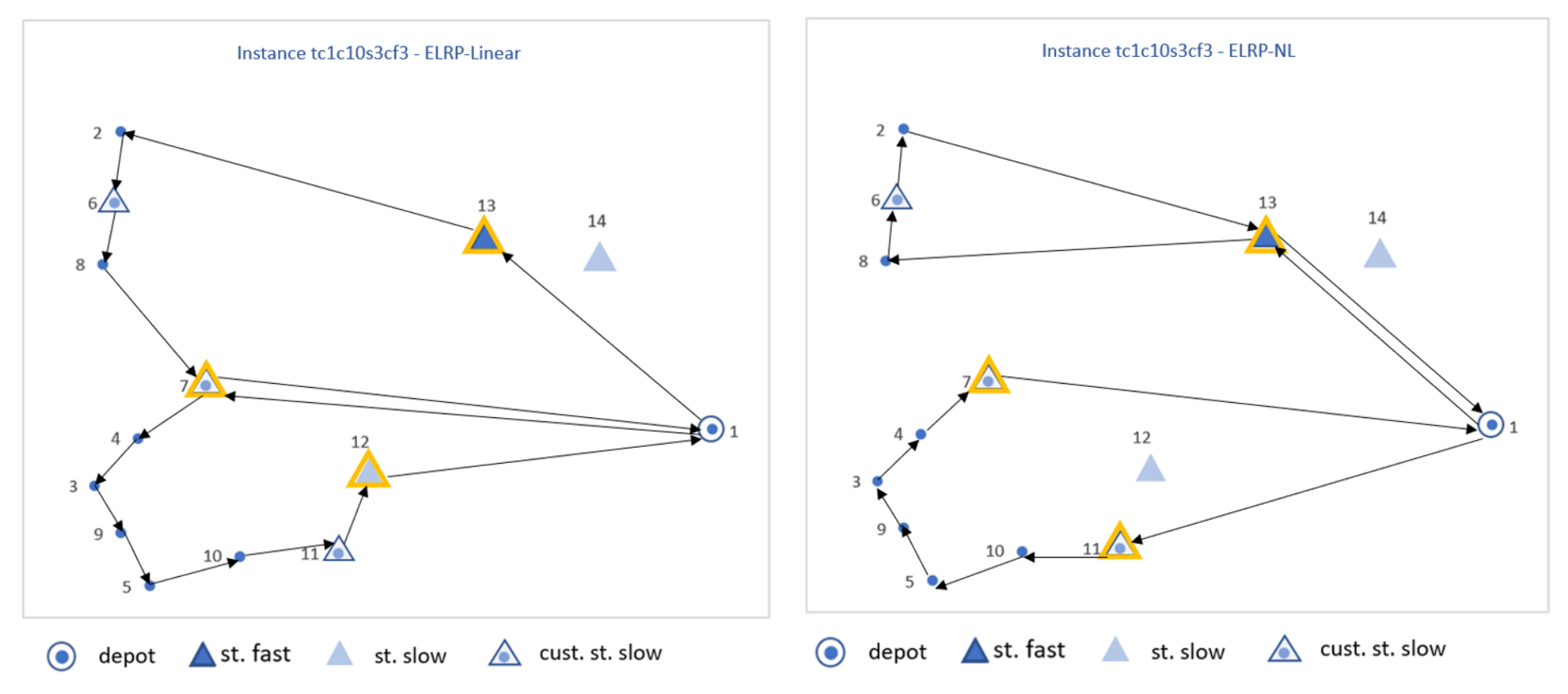}}
\caption{ Detailed solution for instance 12 in ELRP-Linear and ELRP-NLMS\label{fig:5}}
\end{figure}


\section{Conclusions}

In this paper, we introduced electric location-routing models that incorporate nonlinear charging and multiple charging station types to support the planning of EV fleets and charging infrastructure within EV logistic networks. We first extended the classic ELRP formulation proposed by \cite{schiffer2017electric} into an ELRP-NLMS model, and we propose three novel formulations for the problem with alternative approaches for modeling EV battery state-of-charge, time consumption, and nonlinear charging. New instances of the ELRP-NLMS problem have been generated and made available for future studies. We performed extensive computational experiments with 10-80 customer instances, and \change{our analysis demonstrates the effectiveness of the new formulations, reducing the average gap from 29.1\% to 11.9\%, yielding improved solutions for 28 out of 74 instances compared to the node-based formulation}, contributing to the ELRP state-of-the-art formulations available in the literature. The improved formulations optimally solved instances with 10 and 20 customers and found feasible solutions for instances with 40 and 80 customers with a minimum Gap of 9.5\% and 16.4\%, respectively.

The present study also investigates the implications of nonlinear charging in ELRP decision-making. We defined an ELRP model replacing the nonlinear charging with a linear charging function while keeping the same ELRP-NLMS parameters. The results demonstrate that including the nonlinear charging function in ELRP may significantly influence strategic and operational decisions and impact the \change{objective function result}. Therefore, although most ELRP models in the literature assume a linear charging process, considering the nonlinear charging, in addition to being a more realistic approach, is fundamental for more consistent strategic decisions in the optimal design of EV logistic networks.

In future work, developing new metaheuristic solution methods to solve larger instances of the ELRP-NLMS would be fruitful. The improved formulations introduced in this study may also support the development of new exact methods and approaches for ELRP modeling. It would be interesting to consider stochastic elements in the ELRP, such as the uncertainty of the EV range, and investigate how this could impact ELRP strategic decisions. Finally, applying variations of the objective function, such as minimizing the total costs of the network and deriving new extensions that deal with heterogeneous EV fleets and environmental impact, could also lead to new insights into the problem.

\section*{Acknowledgments}
This work was partially supported by CAPES (Coordenacão de Aperfeiçoamento de Pessoal de Níıvel Superior - Finance Code 001), CNPq (Conselho Nacional de Desenvolvimento Científico e Tecnológico - projects 315361/2020-4, 422470/2021-0 and 314420/2023-1), and by FAPERJ (Fundação Carlos Chagas Filho de Amparo à Pesquisa do Estado do Rio de Janeiro - projects E-26/201.417/2022, E-26/010.002232/2019, and E-26/210.041/2023).


\bibliographystyle{unsrt}  
\bibliography{references}  

\clearpage
\appendix
\section{Description of Preprocessing of Recharge Paths}\label{App:A}

In the context of the ELRP-NLMS, we propose a preprocessing procedure to identify strictly worse recharge paths that may be excluded from candidate recharge paths without compromising an optimal solution. We use as a reference to the preprocessing procedure the dominance analysis in recharge paths presented by \cite{froger2019improved} from which we adapted the following definition: A recharge path $p'$ is said to be dominated by another path $p$ between nodes $i,j \in C \cup \{o\}$ if, for every possible SoC at the destination $j$, it is possible for an EV to travel from $i$ to $j$ via path $p$ in a shorter time duration than that of path $p'$ (time-dominance condition), and for every possible time at the destination $j$, it is also possible for an EV to travel from $i$ to $j$ using path $p$ with lower battery consumption compared to path $p'$, reaching $j$ with a higher state of charge (energy-dominance condition). Therefore, if time-dominance and energy-dominance condition is established comparing the paths $p$ and $p'$, for any initial SoC $q_i$ and initial time $\tau_i$, path $p'$ will be strictly worse than path $p$ regarding time and energy performance.

In this Section we detail the two base cases utilized in the Algorithm presented in Section 3.3.1 to identify the time and energy dominance conditions between recharge paths.

\subsection{Base Case 1}

    Consider two paths between nodes $i,j$, where path $p$ has one CS $s$ of type $sp$ and path $p'$ has one CS $s'$ of type $sp'$. The time and energy costs of the first arc of the paths are denoted by the parameters $t_{is}$, $t_{is'}$, $e_{is}$ and $e_{is'}$ respectively. Similarly, the costs for the last arc are $t_{sj}$ , $t_{s'j}$, $e_{sj}$ and $e_{s'j}$. Figure 6 illustrates the graph structure of the two paths.\\

\begin{figure}[htbp]
\centering
    \includegraphics[trim=0 150 0 0, clip, width=0.5\textwidth]{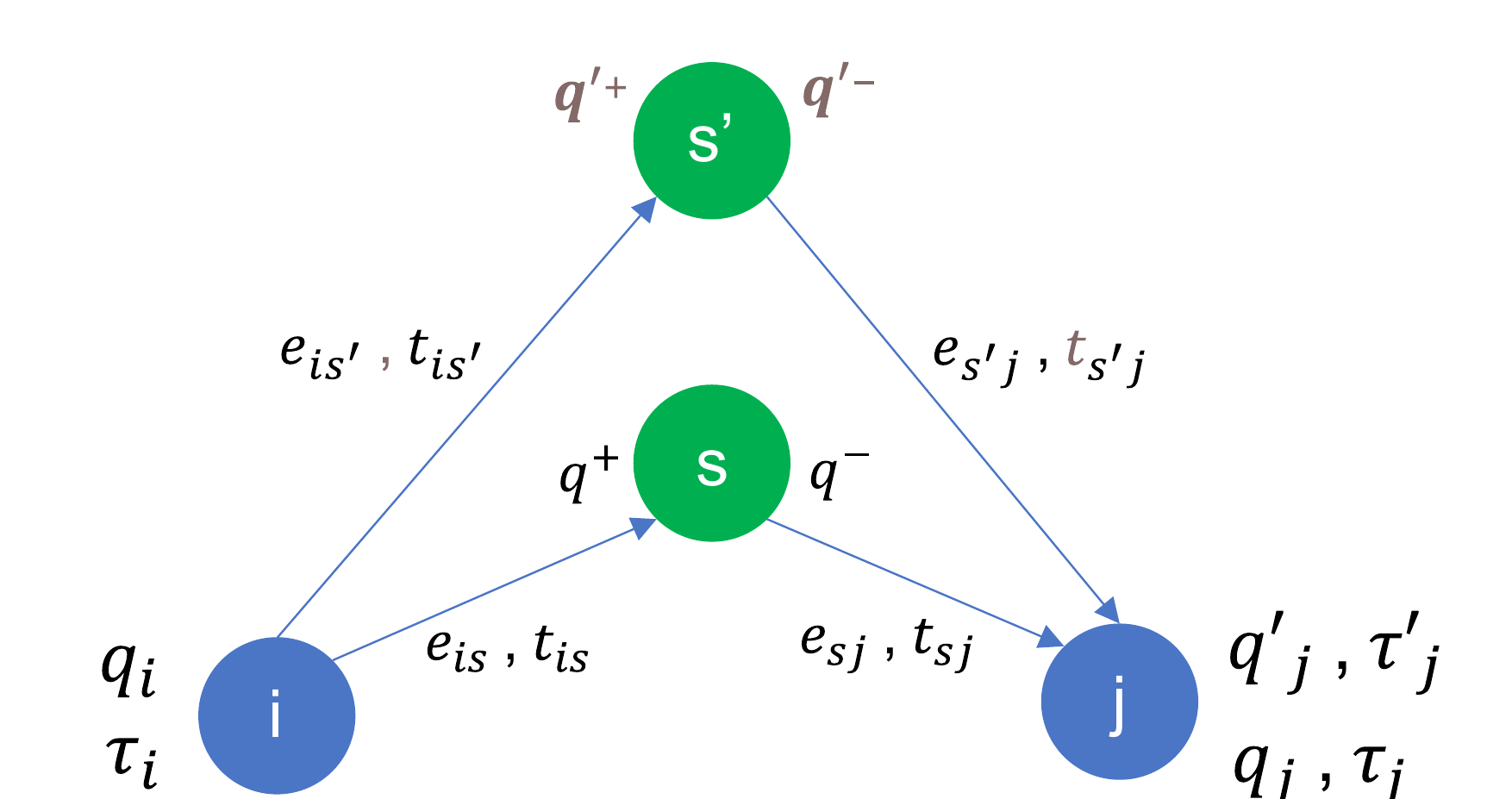}
    \caption{Graph illustration of Base Case 1.}
\label{fig:appendixA}
\end{figure}
    
    \textbf{Case 1 Statement:} If $e_{is} < e_{is'}$ , $e_{sj} < e_{s'j}$ , $(t_{is} + t_{sj}) < (t_{is'} + t_{s'j})$  and the charging speed of CS $s$ is equal to or faster than that of $s'$, then for every possible SoC at the destination $j$, it is possible for an EV to travel from $i$ to $j$ via path $p$ in a shorter time duration than that of path $p'$, and for every possible time at the destination $j$, it is also possible for an EV to travel from $i$ to $j$ using path $p$ with lower battery consumption compared to path $p'$.

    We observe that $e_{is} < e_{is'}$ and $e_{sj} < e_{s'j}$ imply that $t_{is} < t_{is'}$ and $t_{sj} < t_{s'j}$ due to the properties of ELRP-NLMS instance, where the EV speed and the energy consumption rates are assumed to be constant and the arc distances respect the triangular inequality.
    
\textbf{Time-dominance condition in Case 1}

    Consider Case 1 description where an EV with a SoC $q_i$ at node $i$ can travel to node $j$ using the two paths $p$ and $p'$, reaching $j$ with the same SoC as defined by Equation (128).
    The SoC $q^+_s$ and $q^+_{s'}$ of the EV upon entering $s$ and $s'$ are computed by Equations (129) and (130) respectively. Given the condition $e_{is} < e_{is'}$, we can conclude that $q^+_s$ is greater than $q^+_{s'}$ (131). 
    The SoC of the EV when leaving CS $s$ and {s'} is denoted by $q^-_s$ and $q^-_{s'}$ and the SoC of the EV when reaching node $j$ is described by Equations (132) and (133) respectively. Given the condition $e_{sj} < e_{s'j}$ of Case 1 statement, it can be concluded that $q^-_s$ is lower than $q^-_{s'}$ (134).
    
    \vspace*{-15pt}
    \begin{align}
    &q_j = q'_j\label{eq1}\\
    &q^+_s = q_i - e_{is}\label{eq2}\\
    &q^+_{s'} = q_i - e_{is'}\label{eq3}\\
    &q^+_s > q^+_{s'} \label{eq4}\\
    &q_j = q^-_s - e_{sj} \label{eq5}\\
    &q_j = q^-_{s'} - e_{s'j} \label{eq6}\\
    &q^-_s < q^-_{s'} \label{eq7}
    \end{align}

    The SoC of the EV when leaving CS $s$ and $s'$ depends on the amount of energy charged $\Delta q_s$ and $\Delta q_{s'}$  represented by Equations (135) and (136), respectively.
    Based on (131), (134), (135) and (136) we conclude that the amount of energy charged $\Delta q_s$ is lower than $\Delta q_{s'}$, which implies that charging time $\Delta t_s$ is also lower than charging time $\Delta t_{s'}$, given the condition that charging speed of CS $s$ is equal or faster than that of $s'$ of Case 1. This relation is represented by Inequation (137) and can be deduced by the analysis of the charging behavior of the piecewise charging function, illustrated in Figure 7.

    \begin{figure}[htbp]
    \centering
    \includegraphics[width=0.5\textwidth]{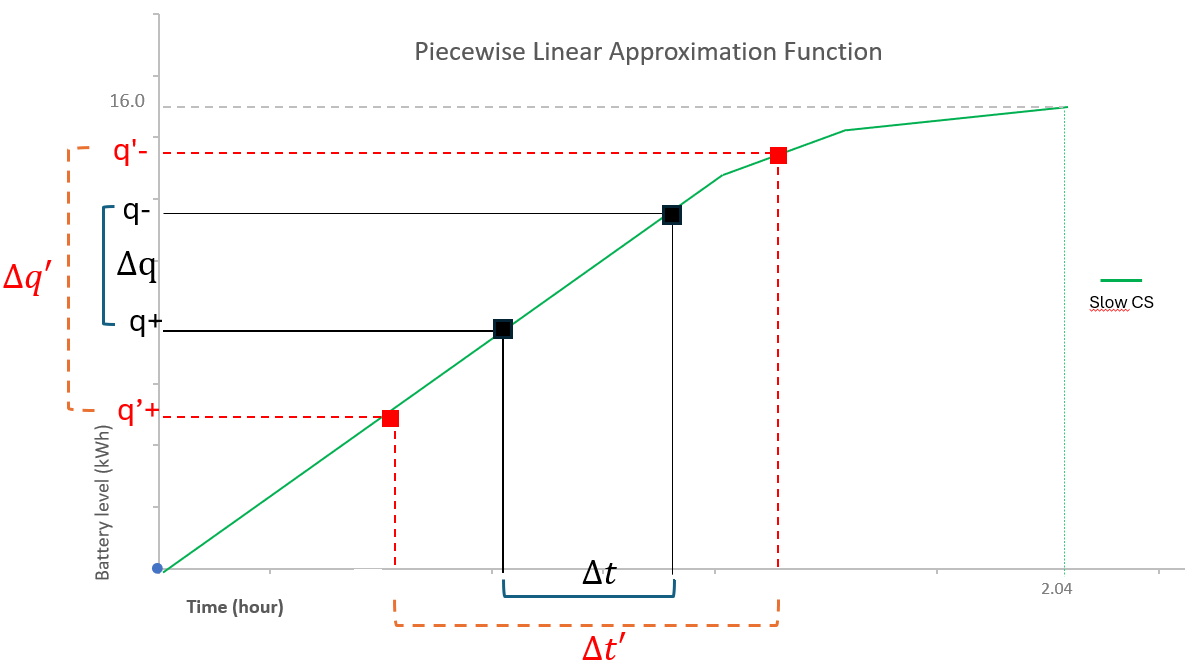}
    \caption{Piecewise linear approximation simulating Case 1 .}
    \label{fig:exemplo}
    \end{figure}
    
   The time in destination $j$ by traveling path $p$ and $p'$ are denoted by $\tau_j$ and $\tau'_j$ respectively, and is calculated by the initial time $\tau_i$ at node $i$ plus the sum of time consumption in each path as described by Equations (138) and (139). Based on (137), (138), (139), and given the condition $(t_{is} + t_{sj}) < (t_{is'} + t_{s'j})$ of Case 1 it can be concluded that the time $\tau_j$ is lower than $\tau'_j$, as described by (140). Therefore, we found time-dominance condition is established in Case 1.
    
    \vspace*{-15pt}
    \begin{align}
    \Delta q_s &= q^-_s - q^+_s \label{eq8}\\
    \Delta q_s' &= q^-_{s'} - q^+_{s'} \label{eq9}\\
    \Delta q_s < \Delta q_{s'} &\implies \Delta t_s < \Delta t_{s'} \label{eq10}
    \end{align}
    \begin{align}
    \tau_j &= \tau_i + t_{is} + t_{sj} + \Delta t_s \label{eq12}\\
    \tau'_j &= \tau_i + t_{is'} + t_{s'j} + \Delta t_s' \label{eq13}\\
    t_{is} + t_{sj} &< t_{is'} + t_{s'j} \implies \tau_j < \tau'_j \label{eq14}
    \end{align}

\textbf{Energy-dominance in Case 1}

    Consider the same context of Case 1, where an EV with a SoC $q_i$ at node $i$ can travel between the two paths $p$ and $p'$, but now reaching destination $j$ at the same time, as defined by (141). Therefore, the time consumption on each path must be equal, as represented by Equation (142). Given the Case 1 condition $(t_{is} + t_{sj}) < (t_{is'} + t_{s'j})$, it can be concluded that $\Delta t_{s}$ must be greater than $\Delta t_{s'}$ to make (142) feasible. This implies that $\Delta q_{s}$ must be greater than $\Delta q_{s'}$ based on Equations (131), (134) and the assumption that the charging speed of CS $s$ is equal to or faster than that of $s'$ in Case 1. Equation (143) describes this relation, which can also be deduced by the analysis of the charging behavior, illustrated in Figure 2.

    The SoC of the EV when reaching node $j$ by paths $p$ and $p'$ are denoted by $q_j$ and $q'_j$ respectively, and is calculated by the initial SoC at node $i$ minus the sum of energy consumption in each path as described by Equations (144) and (145) respectively. Based on (143), (144), (145) and the Case 1 conditions that $e_{is} < e_{is'}$ and $e_{sj} < e_{s'j}$, it can be concluded that the SoC of the EV when reaching destination $q'_j$ is lower than $q_j$. Therefore, we found energy-dominance condition is established in Case 1. Finally, we conclude that under the conditions of Case 1, time-dominance and energy-dominance are established.

\begin{align}
    \tau_j &= \tau'_j\label{eq2}\\
    t_{is} + t_{sj} + \Delta t_s &= t_{is'} + t_{s'j} + \Delta t_{s'}\label{eq15}\\
    \Delta t_s > \Delta t_{s'} &\implies \Delta q_s > \Delta q_{s'} \label{eq16}\\
    q_j &= q_i - e_{is} - e_{sj} + \Delta q_s \label{eq17}\\
    q'_j &= q_i - e'_{is} - e'_{sj} + \Delta q_{s'} \label{eq18}\\
    e_{is} + e_{sj} - \Delta q_s &<  e_{is'} + e_{s'j} - \Delta q_s'\implies q_j > q'_j \label{eq19}
\end{align}

\subsection{Base Case 2}

    Case 2 considers the condition where the recharge speed of CS $s$ in path $p$ is slower than $s'$ in $p$. In this context, Equation (137) may not be universally applicable, and it is necessary to account for the differences between $\Delta t$ and $\Delta q$ when comparing the two paths.
    By analyzing the piecewise charging function, and the graph structure of paths $p$ and $p'$, we can calculate the upper bounds $\bar{\Delta t}$ and $\bar{\Delta q}$, represented in Equations (147) and (148) as the maximum differences in charging time and charged energy between CSs $s$ and $s'$ that an EV may experience when traveling along the two paths. These upper bounds help determine the conditions under which time and energy-dominance between paths can be asserted in Case 2.
    
    \vspace*{-10pt}
    \begin{align}
    \bar{\Delta t} &= max(\Delta t_s - \Delta t_{s'})  & \forall q_j = q'_j \label{eq20}\\
    \bar{\Delta q} &= max(\Delta q_s - \Delta q_{s'})  & \forall \tau_j = \tau'_j 
    \label{eq21}
    \end{align}
    
    Consider two paths between nodes $i,j$, where path $p$ has one CS $s$ of type $sp$, path $p'$ has one CS $s'$ of type $sp'$, and the charging speed of $sp$ is slower than $sp'$. Consider $\bar{\Delta t}$ and $\bar{\Delta q}$ as the maximum difference in charging time and charged energy between $sp$ and $sp'$. Figure 8 illustrates the graph structure of the two paths.
    
    \begin{figure}[htbp]
        \centering
        \includegraphics[trim=0 150 0 0, clip, width=0.5\textwidth]{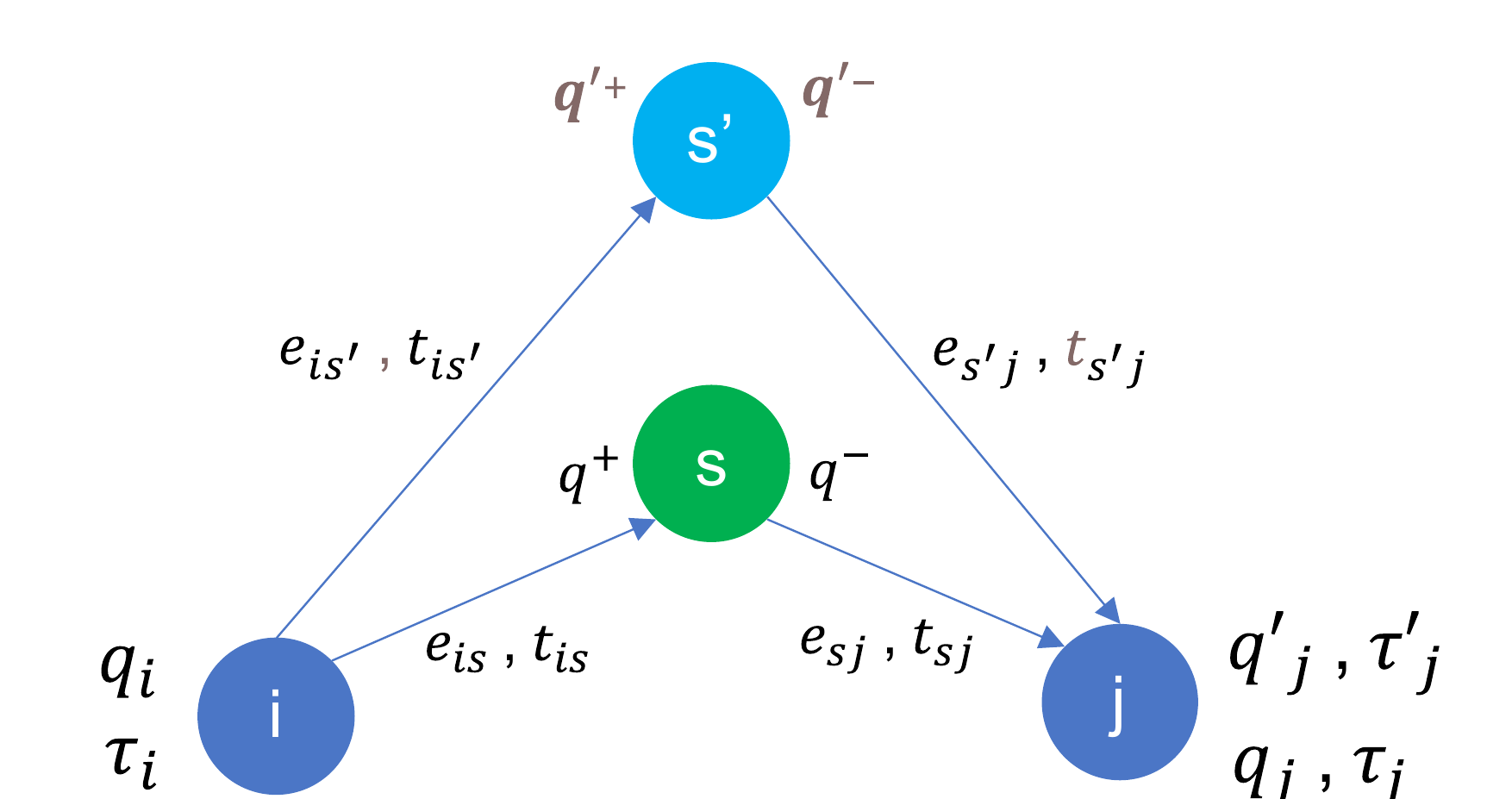}
        \caption{Graph Structure - Case 2}
        \label{fig:exemplo}
    \end{figure}
    
    \textbf{Case 2 Statement:} In the case that the charging speed of CS $s$ is slower than that of $s'$, if $e_{is} < e_{is'}$ , $e_{sj} < e_{s'j}$ , $t_{is} + t_{sj} + \bar{\Delta t} < t_{is'} + t_{s'j}$ , $e_{is} + e_{sj} + \bar{\Delta q} < e_{is'} + e_{s'j}$, then for every possible SoC at the destination $j$, it is possible for an EV to travel from $i$ to $j$ via path $p$ in a shorter time duration than that of path $p'$, and for every possible time at the destination $j$, it is also possible for an EV to travel from $i$ to $j$ using path $p$ with lower battery consumption compared to path $p'$.

    \textbf{Time dominance in Case 2}

    Referencing (128)-(136), (138) and (139) of Base Case 1, consider Case 2 description where an EV with a SoC $q_i$ at node $i$ can travel to node $j$ using the two paths $p$ and $p'$, reaching $j$ with the same SoC $q_j = q'_j$.

    The $\bar{\Delta t}$ can be calculated by comparing the two piecewise functions and finding the highest time difference between each function considering both paths. Based on path structure, it is possible to identify the bounds of $q^+$ and $q^-$ as defined in (129), (130), and identify the maximum $\Delta t$ among all possible charging amount as illustrated by Figure 9. Considering the difference in $q^+$ as defined in (131), the $\Delta t$ can be straightened, by shifting the piecewise function of the slower CS $s$ by the equivalent charging time of $q^+$, $\Delta t_{in}$, as illustrated in Figure 10.
    
    \begin{figure}[htbp]
        \centering
        \includegraphics[width=0.5\textwidth]{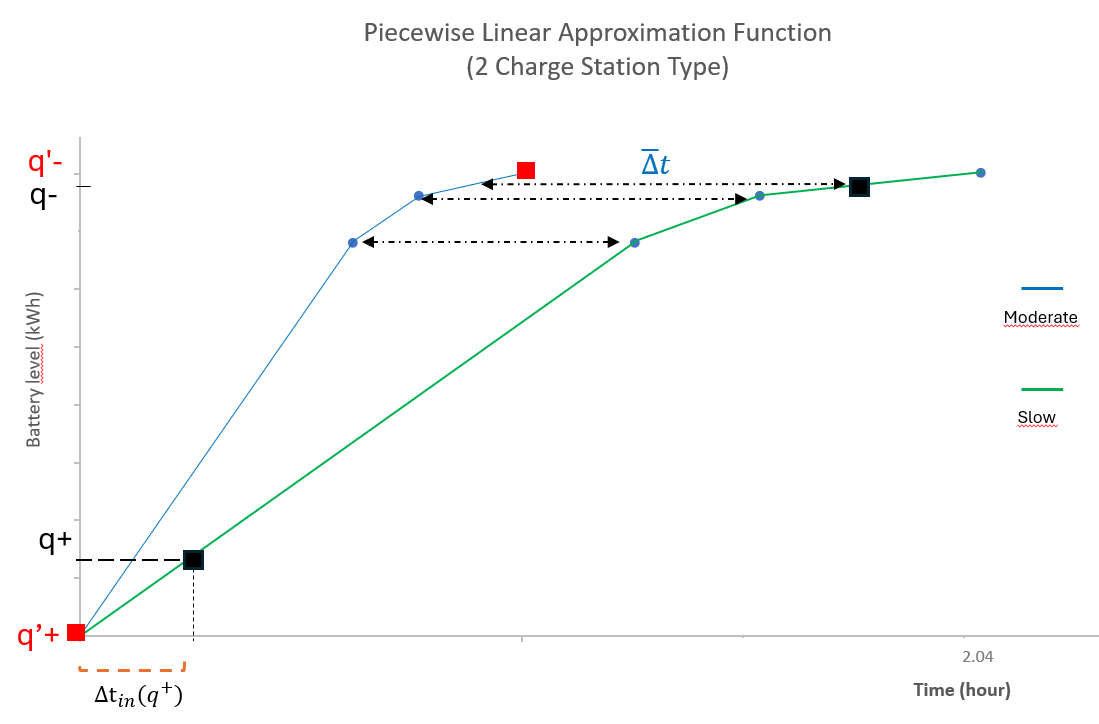}
        \caption{Calculating $\bar{\Delta t}$ in Case 2}
        \label{fig:exemplo}
    \end{figure}

    \begin{figure}[htbp]
        \centering
        \includegraphics[width=0.5\textwidth]{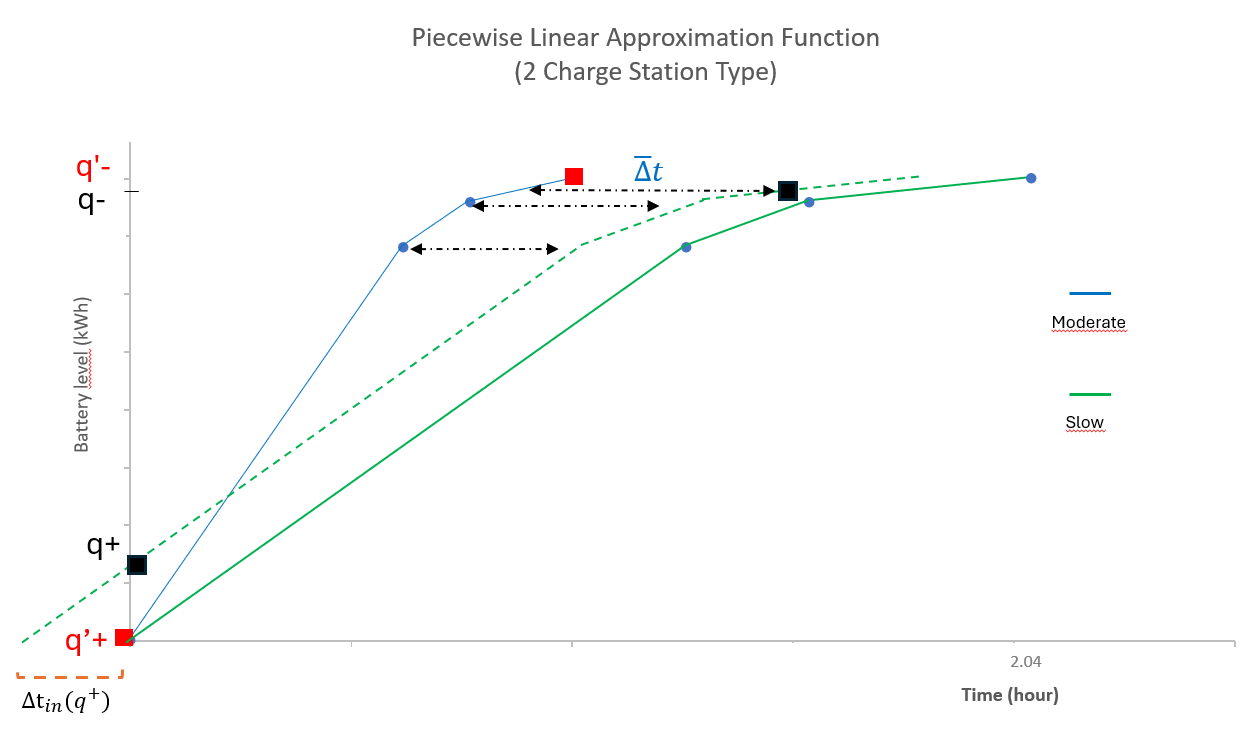}
        \caption{Straightening $\bar{\Delta t}$ in Case 2}
        \label{fig:exemplo}
    \end{figure}
    
    Based on the definition of $\bar{\Delta t}$ as the maximum difference in charging time between $s$ and $s'$ in path $p$ and $p'$ respectively, we can assert that ${\Delta t_s}$ will necessarily be lower or equal than ${\Delta t_s'} + \bar{\Delta t}$ as defined in (149). Considering time consumption Equations (11) and (12), Equation (150) shows that if the time consumption in path $p$ is lower than $p'$, this implies in $\tau_j < \tau'_j$. In Equation (151) we substitute the $\bar{\Delta t}$ by ${\Delta t_s'} + \bar{\Delta t}$. Finally, Equation (152) aligns with the statement of Case 2 statement, founding that time-dominance condition is established in Case 2.

    \begin{align}
    \Delta t_s &\leq \Delta t_{s'} + \bar{\Delta t} \label{eq22}\\
    t_{is} + t_{sj} + \Delta t_s &< t_{is'} + t_{s'j} + \Delta t_s' \implies \tau_j < \tau'_j \label{eq23}\\
    t_{is} + t_{sj} + \Delta t_s' + \bar{\Delta t} &< t_{is'} + t_{s'j} + \Delta t_s' \implies \tau_j < \tau'_j \label{eq24}\\
    t_{is} + t_{sj} + \bar{\Delta t} &< t_{is'} + t_{s'j} \implies \tau_j < \tau'_j \label{eq25}
    \end{align}

    \textbf{Energy dominance in Case 2}

    Consider the same context of Case 2, where an EV with a SoC $q_i$ at node $i$ can travel to node $j$ using the two paths $p$ and $p'$, but now reaching destination $j$ at the same time $\tau_j = \tau'_j$.

    Similar to the calculation of $\bar{\Delta t}$, the upper bound $\bar{\Delta q}$ can be calculated by comparing the two piecewise functions and finding the highest charge amount difference between each function considering both paths. Based on the path structure, it is possible to identify the maximum $\Delta q$ among all possible charging amounts, as illustrated in Figure 11. The $\bar{\Delta q}$ can also be straightened, by shifting the piecewise function of the slower CS $s$ by the equivalent charging time of $q^+$, as illustrated in Figure 12.

    \begin{figure}[htbp]
        \centering
        \includegraphics[width=0.5\textwidth]{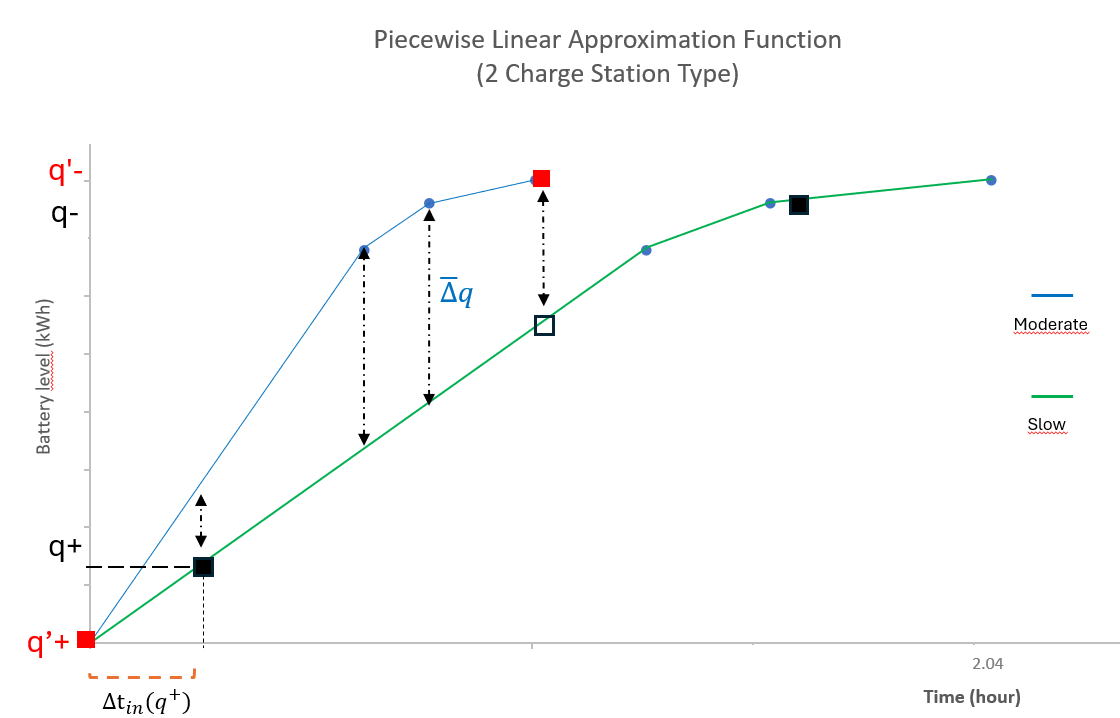}
        \caption{Calculating $\bar{\Delta q}$ in Case 2}
        \label{fig:exemplo}
    \end{figure}

    \begin{figure}[htbp]
        \centering
        \includegraphics[width=0.5\textwidth]{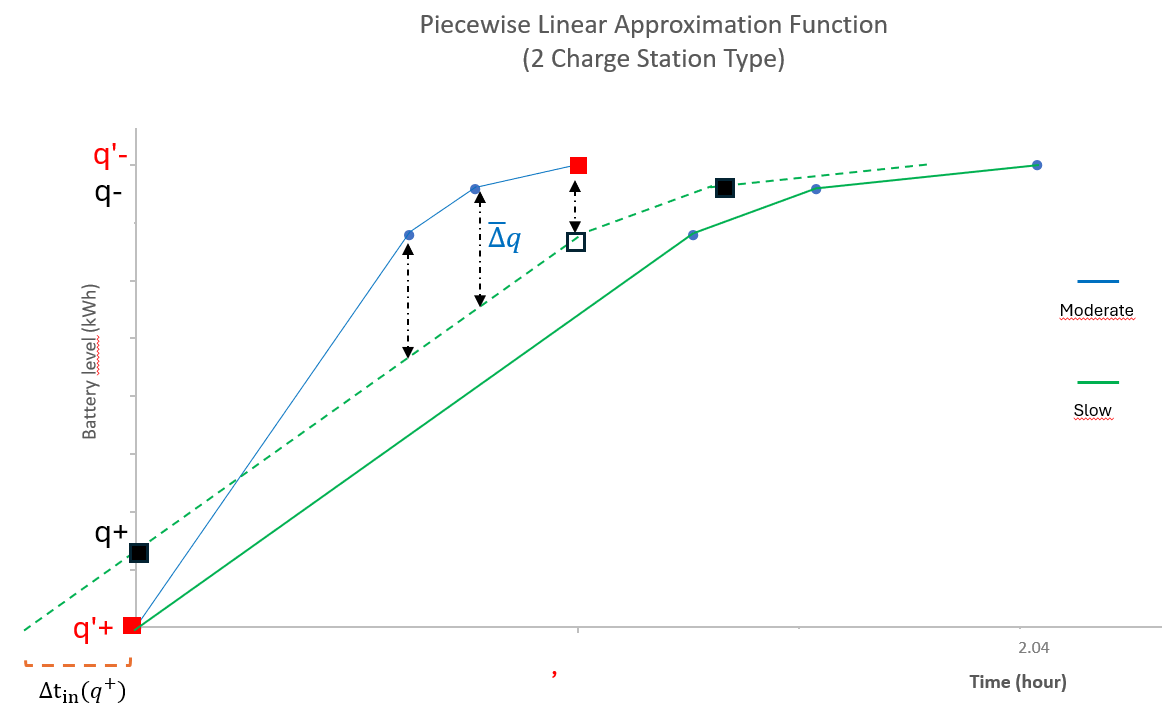}
        \caption{Straightening $\bar{\Delta q}$ in Case 2}
        \label{fig:exemplo}
    \end{figure}
    
    Based on the definition of $\bar{\Delta q}$ as the maximum difference in charge amount between $s$ and $s'$ in path $p$ and $p'$ respectively, we can assert that ${\Delta q_s}$ will necessarily be lower or equal than ${\Delta q_s'} + \bar{\Delta q}$ as defined in (153).  Considering energy consumption Equations (144) and (145), Equation (154) shows that if the energy consumption in path $p$ is higher than $p'$, this implies in the SoC $q_j$ will be higher than $q'_j$. In (155) we substitute $\bar{\Delta q}$ by ${\Delta q'_s} + \bar{\Delta q}$. Finally, Equation (156) aligns with the statement of Case 2, founding that energy-dominance condition is established in Case 2.

    \begin{align}
    \Delta q_s &\leq  \Delta q_{s'} + \bar{\Delta q} \label{eq26}\\
    e_{is} + e_{sj} - \Delta q_s &< e_{is'} + e_{s'j} - \Delta q_s' \implies q_j > q'_j \label{eq27}\\
    e_{is} + e_{sj} + \Delta q_s' + \bar{\Delta q} &< e_{is'} + e_{s'j} + \Delta q_s' \implies q_j > q'_j \label{eq28}\\
    e_{is} + e_{sj} + \bar{\Delta q} &< e_{is'} + e_{s'j} \implies q_j > q'_j \label{eq29}
    \end{align}

\clearpage
\section{Detailed Results}\label{App:B}

This section presents the detailed results for each of the 80 instances and the four different models, computing the objective function value, the MIP Gap, and the runtime. TL states that the model instance reached the time limit of 10800 seconds.

\begin{table}[htbp]
  \centering
    \caption{Results for the 10 customer instances.\label{tab:A1}}
  \scriptsize
    \begin{tabular}{cccrrccrcccrcccrcc}
    \toprule
    \multicolumn{2}{p{7.57em}}{\textbf{10-customer}} &       & \multicolumn{3}{p{11.14em}}{\textbf{M1}} &       & \multicolumn{3}{p{9.285em}}{\textbf{M2}} &       & \multicolumn{3}{p{8.78em}}{\textbf{M3}} &       & \multicolumn{3}{p{8.43em}}{\textbf{M4}} \\
\cmidrule{1-2}\cmidrule{4-6}\cmidrule{8-10}\cmidrule{12-14}\cmidrule{16-18}    \textbf{\#} & \textbf{Instances} &       & \multicolumn{1}{c}{\textbf{Obj}} & \multicolumn{1}{c}{\textbf{Gap}} & \textbf{Time} &       & \multicolumn{1}{c}{\textbf{Obj}} & \textbf{Gap} & \textbf{Time} &       & \multicolumn{1}{c}{\textbf{Obj}} & \textbf{Gap} & \textbf{Time} &       & \multicolumn{1}{c}{\textbf{Obj}} & \textbf{Gap} & \textbf{Time} \\
\cmidrule{1-2}\cmidrule{4-6}\cmidrule{8-10}\cmidrule{12-14}\cmidrule{16-18}    1     & tc0c10s2cf1-p4 &       & 19.753 & 0.0\% & 94    &       & 19.753 & 0.0\% & 784   &       & 19.753 & 0.0\% & 84    &       & 19.753 & 0.0\% & 109 \\
    2     & tc0c10s2ct1-p4 &       & 11.378 & 0.0\% & 210   &       & 11.378 & 0.0\% & 70    &       & 11.378 & 0.0\% & 72    &       & 11.378 & 0.0\% & 149 \\
    3     & tc0c10s3cf1-p6 &       & 10.921 & 0.0\% & 57    &       & 10.921 & 0.0\% & 66    &       & 10.921 & 0.0\% & 32    &       & 10.921 & 0.0\% & 57 \\
    4     & tc0c10s3ct1-p6 &       & 10.536 & 0.0\% & 26    &       & 10.536 & 0.0\% & 19    &       & 10.536 & 0.0\% & 26    &       & 10.536 & 0.0\% & 48 \\
    5     & tc1c10s2cf2-p4 &       & 9.034 & 0.0\% & 1062  &       & 9.034 & 0.0\% & 29    &       & 9.034 & 0.0\% & 87    &       & 9.034 & 0.0\% & 133 \\
    6     & tc1c10s2cf3-p4 &       & 12.583 & 0.0\% & 438   &       & 12.583 & 0.0\% & 15    &       & 12.583 & 0.0\% & 11    &       & 12.583 & 0.0\% & 17 \\
    7     & tc1c10s2cf4-p4 &       & 16.097 & 0.0\% & 6131  &       & 16.097 & 0.0\% & 145   &       & 16.097 & 0.0\% & 87    &       & 16.097 & 0.0\% & 180 \\
    8     & tc1c10s2ct2-p4 &       & 10.466 & 24.4\% & TL    &       & 9.199 & 0.0\% & 189   &       & 9.199 & 0.0\% & 780   &       & 9.199 & 0.0\% & 367 \\
    9     & tc1c10s2ct3-p4 &       & 12.307 & 24.0\% & TL    &       & 12.307 & 0.0\% & 127   &       & 12.307 & 0.0\% & 27    &       & 12.307 & 0.0\% & 65 \\
    10    & tc1c10s2ct4-p4 &       & 13.826 & 0.0\% & 114   &       & 13.826 & 0.0\% & 50    &       & 13.826 & 0.0\% & 64    &       & 13.826 & 0.0\% & 65 \\
    11    & tc1c10s3cf2-p6 &       & 9.034 & 0.0\% & 5370  &       & 9.034 & 0.0\% & 22    &       & 9.034 & 0.0\% & 129   &       & 9.034 & 0.0\% & 1327 \\
    12    & tc1c10s3cf3-p6 &       & 12.583 & 27.8\% & TL    &       & 12.583 & 0.0\% & 46    &       & 12.583 & 0.0\% & 49    &       & 12.583 & 0.0\% & 28 \\
    13    & tc1c10s3cf4-p6 &       & 14.902 & 0.0\% & 3186  &       & 14.902 & 0.0\% & 121   &       & 14.902 & 0.0\% & 82    &       & 14.902 & 0.0\% & 277 \\
    14    & tc1c10s3ct2-p6 &       & 9.199 & 11.3\% & TL    &       & 9.199 & 0.0\% & 88    &       & 9.199 & 0.0\% & 207   &       & 9.199 & 0.0\% & 1625 \\
    15    & tc1c10s3ct3-p6 &       & 12.932 & 20.2\% & TL    &       & 12.932 & 0.0\% & 33    &       & 12.932 & 0.0\% & 129   &       & 12.932 & 0.0\% & 335 \\
    16    & tc1c10s3ct4-p6 &       & 13.205 & 0.0\% & 85    &       & 13.205 & 0.0\% & 38    &       & 13.205 & 0.0\% & 126   &       & 13.205 & 0.0\% & 299 \\
    17    & tc2c10s2cf0-p4 &       & 11.345 & 33.2\% & TL    &       & 11.345 & 0.0\% & 514   &       & 11.345 & 0.0\% & 355   &       & 11.345 & 0.0\% & 1054 \\
    18    & tc2c10s2ct0-p4 &       & 11.345 & 34.8\% & TL    &       & 11.345 & 0.0\% & 5940  &       & 11.345 & 0.0\% & 1019  &       & 11.345 & 0.0\% & 2445 \\
    19    & tc2c10s3cf0-p6 &       & 11.285 & 34.2\% & TL    &       & 11.285 & 0.0\% & 514   &       & 11.285 & 0.0\% & 837   &       & 11.285 & 0.0\% & 5867 \\
    20    & tc2c10s3ct0-p6 &       & 11.285 & 34.3\% & TL    &       & 11.285 & 0.0\% & 65    &       & 11.285 & 0.0\% & 2751  &       & 11.285 & 0.0\% & 1129 \\
    \bottomrule
    \end{tabular}%
  \label{tab:addlabel}%
\end{table}%

\begin{table}[htbp]
  \centering
  \caption{Results for the 20 customer instances.\label{tab:A2}}
  \scriptsize
    \begin{tabular}{cccrrccrrccrrccrrc}
    \toprule
    \multicolumn{2}{p{7.57em}}{\textbf{20-customer}} &       & \multicolumn{3}{p{11.14em}}{\textbf{M1}} &       & \multicolumn{3}{p{9.285em}}{\textbf{M2}} &       & \multicolumn{3}{p{8.78em}}{\textbf{M3}} &       & \multicolumn{3}{p{8.43em}}{\textbf{M4}} \\
\cmidrule{1-2}\cmidrule{4-6}\cmidrule{8-10}\cmidrule{12-14}\cmidrule{16-18}    \textbf{\#} & \textbf{Instances} &       & \multicolumn{1}{c}{\textbf{Obj}} & \multicolumn{1}{c}{\textbf{Gap}} & \textbf{Time} &       & \multicolumn{1}{c}{\textbf{Obj}} & \multicolumn{1}{c}{\textbf{Gap}} & \textbf{Time} &       & \multicolumn{1}{c}{\textbf{Obj}} & \multicolumn{1}{c}{\textbf{Gap}} & \textbf{Time} &       & \multicolumn{1}{c}{\textbf{Obj}} & \multicolumn{1}{c}{\textbf{Gap}} & \textbf{Time} \\
\cmidrule{1-2}\cmidrule{4-6}\cmidrule{8-10}\cmidrule{12-14}\cmidrule{16-18}    21    & tc2c20s3cf0-p6 &       & 24.464 & 62.2\% & TL    &       & 24.372 & 36.4\% & TL    &       & 24.040 & 31.8\% & TL    &       & 24.040 & 32.6\% & TL \\
    22    & tc1c20s3cf1-p6 &       & 14.747 & 27.9\% & TL    &       & 14.747 & 12.0\% & TL    &       & 14.747 & 6.4\% & TL    &       & 14.747 & 9.5\% & TL \\
    23    & tc0c20s3cf2-p6 &       & 16.651 & 20.5\% & TL    &       & 16.671 & 11.5\% & TL    &       & 16.651 & 17.1\% & TL    &       & 16.921 & 20.7\% & TL \\
    24    & tc1c20s3cf3-p6 &       & 11.861 & 20.4\% & TL    &       & 11.861 & 0.0\% & 106   &       & 11.861 & 0.0\% & 92    &       & 11.861 & 0.0\% & 275 \\
    25    & tc1c20s3cf4-p6 &       & 16.095 & 16.0\% & TL    &       & 16.088 & 0.0\% & 2127  &       & 16.088 & 0.0\% & 4085  &       & 16.088 & 0.0\% & 9803 \\
    26    & tc2c20s3ct0-p6 &       & 24.464 & 62.4\% & TL    &       & 24.307 & 35.8\% & TL    &       & 24.275 & 31.1\% & TL    &       & 24.304 & 32.6\% & TL \\
    27    & tc1c20s3ct1-p6 &       & 15.116 & 30.7\% & TL    &       & 15.083 & 11.8\% & TL    &       & 15.190 & 9.3\% & TL    &       & 15.116 & 10.7\% & TL \\
    28    & tc0c20s3ct2-p6 &       & 16.627 & 21.7\% & TL    &       & 16.627 & 11.8\% & TL    &       & 16.627 & 18.1\% & TL    &       & 16.627 & 19.8\% & TL \\
    29    & tc1c20s3ct3-p6 &       & 11.976 & 20.9\% & TL    &       & 11.861 & 0.0\% & 35    &       & 11.861 & 0.0\% & 99    &       & 11.861 & 0.0\% & 192 \\
    30    & tc1c20s3ct4-p6 &       & 15.792 & 15.7\% & TL    &       & 15.792 & 0.0\% & 608   &       & 15.792 & 0.0\% & 2063  &       & 15.792 & 0.0\% & 9982 \\
    31    & tc2c20s4cf0-p8 &       & 25.108 & 62.0\% & TL    &       & 24.741 & 36.0\% & TL    &       & 24.672 & 30.8\% & TL    &       & 24.723 & 33.9\% & TL \\
    32    & tc1c20s4cf1-p8 &       & 16.099 & 20.7\% & TL    &       & 16.041 & 18.3\% & TL    &       & 16.143 & 17.4\% & TL    &       & 16.212 & 22.0\% & TL \\
    33    & tc0c20s4cf2-p8 &       & 16.206 & 27.7\% & TL    &       & 16.206 & 15.6\% & TL    &       & 16.206 & 15.3\% & TL    &       & 16.206 & 17.8\% & TL \\
    34    & tc1c20s4cf3-p8 &       & 13.066 & 26.7\% & TL    &       & 12.967 & 0.0\% & 1717  &       & 12.967 & 0.0\% & 1148  &       & 12.967 & 0.0\% & 5221 \\
    35    & tc1c20s4cf4-p8 &       & 15.825 & 14.1\% & TL    &       & 15.825 & 0.0\% & 2437  &       & 15.825 & 0.0\% & 6432  &       & 15.825 & 5.8\% & TL \\
    36    & tc2c20s4ct0-p8 &       & 26.200 & 58.0\% & TL    &       & 26.185 & 42.9\% & TL    &       & 26.018 & 35.2\% & TL    &       & 26.120 & 38.1\% & TL \\
    37    & tc1c20s4ct1-p8 &       & 16.070 & 26.3\% & TL    &       & 16.070 & 19.5\% & TL    &       & 16.070 & 16.8\% & TL    &       & 16.375 & 22.6\% & TL \\
    38    & tc0c20s4ct2-p8 &       & 16.056 & 29.5\% & TL    &       & 16.056 & 14.4\% & TL    &       & 16.056 & 15.8\% & TL    &       & 16.056 & 17.5\% & TL \\
    39    & tc1c20s4ct3-p8 &       & 12.866 & 27.5\% & TL    &       & 12.866 & 0.0\% & 3266  &       & 12.866 & 0.0\% & 2109  &       & 12.866 & 0.0\% & 7823 \\
    40    & tc1c20s4ct4-p8 &       & 15.836 & 16.4\% & TL    &       & 15.825 & 0.0\% & 2160  &       & 15.825 & 2.1\% & TL    &       & 15.825 & 6.7\% & TL \\
    \bottomrule
    \end{tabular}%
  \label{tab:addlabel}%
\end{table}%

\begin{table}[htbp]
  \centering
    \caption{Results for the 40 customer instances.\label{tab:A3}}
  \scriptsize
    \begin{tabular}{cccrrccrrccrrccrrc}
    \toprule
    \multicolumn{2}{p{7.57em}}{\textbf{40-customer}} &       & \multicolumn{3}{p{11.14em}}{\textbf{M1}} &       & \multicolumn{3}{p{9.285em}}{\textbf{M2}} &       & \multicolumn{3}{p{8.78em}}{\textbf{M3}} &       & \multicolumn{3}{p{8.43em}}{\textbf{M4}} \\
\cmidrule{1-2}\cmidrule{4-6}\cmidrule{8-10}\cmidrule{12-14}\cmidrule{16-18}    \textbf{\#} & \textbf{Instances} &       & \multicolumn{1}{c}{\textbf{Obj}} & \multicolumn{1}{c}{\textbf{Gap}} & \textbf{Time} &       & \multicolumn{1}{c}{\textbf{Obj}} & \multicolumn{1}{c}{\textbf{Gap}} & \textbf{Time} &       & \multicolumn{1}{c}{\textbf{Obj}} & \multicolumn{1}{c}{\textbf{Gap}} & \textbf{Time} &       & \multicolumn{1}{c}{\textbf{Obj}} & \multicolumn{1}{c}{\textbf{Gap}} & \textbf{Time} \\
\cmidrule{1-2}\cmidrule{4-6}\cmidrule{8-10}\cmidrule{12-14}\cmidrule{16-18}    41    & tc2c40s5cf2-p10 &       & 27.425 & \multicolumn{1}{c}{55.2\%} & TL    &       & \multicolumn{1}{c}{-} & \multicolumn{1}{c}{-} & TL    &       & \multicolumn{1}{c}{-} & \multicolumn{1}{c}{-} & TL    &       & 28.187 & 37.8\% & TL \\
    42    & tc2c40s5cf3-p10 &       & 19.171 & 41.9\% & TL    &       & 19.171 & 24.6\% & TL    &       & 19.246 & 22.5\% & TL    &       & 19.991 & 34.4\% & TL \\
    43    & tc0c40s5cf4-p10 &       & 33.034 & 50.8\% & TL    &       & \multicolumn{1}{c}{-} & \multicolumn{1}{c}{-} & TL    &       & \multicolumn{1}{c}{-} & \multicolumn{1}{c}{-} & TL    &       & \multicolumn{1}{c}{-} & \multicolumn{1}{c}{-} & TL \\
    44    & tc2c40s8cf2-p16 &       & 26.781 & \multicolumn{1}{c}{55.3\%} & TL    &       & 26.705 & 26.0\% & TL    &       & 27.034 & 28.9\% & TL    &       & 28.359 & 39.0\% & TL \\
    45    & tc2c40s8cf3-p16 &       & 20.522 & 46.1\% & TL    &       & 19.642 & 25.0\% & TL    &       & 20.207 & 29.4\% & TL    &       & 23.410 & 44.0\% & TL \\
    46    & tc0c40s8ct0-p16 &       & 25.422 & \multicolumn{1}{c}{39.2\%} & TL    &       & 24.802 & 12.0\% & TL    &       & 25.455 & 13.2\% & TL    &       & 31.191 & 37.3\% & TL \\
    47    & tc2c40s8ct3-p16 &       & 22.545 & 50.7\% & TL    &       & 22.571 & 33.5\% & TL    &       & 24.371 & 47.8\% & TL    &       & 29.880 & 61.5\% & TL \\
    48    & tc0c40s5cf0-p10 &       & 29.097 & 44.9\% & TL    &       & 27.245 & 14.9\% & TL    &       & 27.658 & 17.1\% & TL    &       & 28.012 & 20.7\% & TL \\
    49    & tc1c40s5cf1-p10 &       & \multicolumn{1}{c}{-} & \multicolumn{1}{c}{-} & TL    &       & \multicolumn{1}{c}{-} & \multicolumn{1}{c}{-} & TL    &       & \multicolumn{1}{c}{-} & \multicolumn{1}{c}{-} & TL    &       & \multicolumn{1}{c}{-} & \multicolumn{1}{c}{-} & TL \\
    50    & tc0c40s5ct0-p10 &       & 27.027 & 42.2\% & TL    &       & 26.186 & 13.4\% & TL    &       & 26.693 & 14.8\% & TL    &       & 28.613 & 26.8\% & TL \\
    51    & tc1c40s5ct1-p10 &       & \multicolumn{1}{c}{-} & \multicolumn{1}{c}{-} & TL    &       & \multicolumn{1}{c}{-} & \multicolumn{1}{c}{-} & TL    &       & \multicolumn{1}{c}{-} & \multicolumn{1}{c}{-} & TL    &       & \multicolumn{1}{c}{-} & \multicolumn{1}{c}{-} & TL \\
    52    & tc2c40s5ct2-p10 &       & 26.216 & 53.8\% & TL    &       & 26.909 & 30.3\% & TL    &       & 26.905 & 28.4\% & TL    &       & 27.290 & 38.9\% & TL \\
    53    & tc2c40s5ct3-p10 &       & 18.756 & 41.1\% & TL    &       & 18.714 & 28.0\% & TL    &       & 19.093 & 23.1\% & TL    &       & 19.324 & 35.2\% & TL \\
    54    & tc0c40s5ct4-p10 &       & 31.707 & 49.7\% & TL    &       & \multicolumn{1}{c}{-} & \multicolumn{1}{c}{-} & TL    &       & \multicolumn{1}{c}{-} & \multicolumn{1}{c}{-} & TL    &       & \multicolumn{1}{c}{-} & \multicolumn{1}{c}{-} & TL \\
    55    & tc0c40s8cf0-p16 &       & 26.278 & \multicolumn{1}{c}{40.7\%} & TL    &       & 25.930 & 11.0\% & TL    &       & 26.018 & 11.9\% & TL    &       & 27.243 & 23.7\% & TL \\
    56    & tc1c40s8cf1-p16 &       & \multicolumn{1}{c}{-} & \multicolumn{1}{c}{-} & TL    &       & \multicolumn{1}{c}{-} & \multicolumn{1}{c}{-} & TL    &       & \multicolumn{1}{c}{-} & \multicolumn{1}{c}{-} & TL    &       & \multicolumn{1}{c}{-} & \multicolumn{1}{c}{-} & TL \\
    57    & tc0c40s8cf4-p16 &       & 30.303 & 48.4\% & TL    &       & 31.128 & 34.8\% & TL    &       & \multicolumn{1}{c}{28.227} & \multicolumn{1}{c}{0.2701} & TL    &       & \multicolumn{1}{c}{-} & \multicolumn{1}{c}{-} & TL \\
    58    & tc1c40s8ct1-p16 &       & \multicolumn{1}{c}{-} & \multicolumn{1}{c}{-} & TL    &       & \multicolumn{1}{c}{-} & \multicolumn{1}{c}{-} & TL    &       & \multicolumn{1}{c}{-} & \multicolumn{1}{c}{-} & TL    &       & \multicolumn{1}{c}{-} & \multicolumn{1}{c}{-} & TL \\
    59    & tc2c40s8ct2-p16 &       & 26.180 & \multicolumn{1}{c}{54.3\%} & TL    &       & 26.107 & 25.2\% & TL    &       & 26.863 & 34.0\% & TL    &       & 32.945 & 50.1\% & TL \\
    60    & tc0c40s8ct4-p16 &       & 28.819 & 45.7\% & TL    &       & 28.971 & 28.6\% & TL    &       & 28.502 & 31.8\% & TL    &       & 32.818 & 42.5\% & TL \\
    \bottomrule
    \end{tabular}%
  \label{tab:addlabel}%
\end{table}%

\begin{table}[htbp]
  \centering
    \caption{Results for the 80 customer instances.\label{tab:A4}}
  \scriptsize
    \begin{tabular}{cccrrccccccccccccc}
    \toprule
    \multicolumn{2}{p{7.57em}}{\textbf{80-customer}} &       & \multicolumn{3}{p{11.14em}}{\textbf{M1}} &       & \multicolumn{3}{p{9.285em}}{\textbf{M2}} &       & \multicolumn{3}{p{8.78em}}{\textbf{M3}} &       & \multicolumn{3}{p{8.43em}}{\textbf{M4}} \\
\cmidrule{1-2}\cmidrule{4-6}\cmidrule{8-10}\cmidrule{12-14}\cmidrule{16-18}    \textbf{\#} & \textbf{Instances} &       & \multicolumn{1}{c}{\textbf{Obj}} & \multicolumn{1}{c}{\textbf{Gap}} & \textbf{Time} &       & \textbf{Obj} & \textbf{Gap} & \textbf{Time} &       & \textbf{Obj} & \textbf{Gap} & \textbf{Time} &       & \textbf{Obj} & \textbf{Gap} & \textbf{Time} \\
\cmidrule{1-2}\cmidrule{4-6}\cmidrule{8-10}\cmidrule{12-14}\cmidrule{16-18}    61    & tc0c80s8cf0-p16 &       & 38.205 & 45.3\% & TL    &       & -     & -     & TL    &       & -     & -     & TL    &       & -     & -     & TL \\
    62    & tc0c80s8cf1-p16 &       & 55.182 & 60.5\% & TL    &       & -     & -     & TL    &       & -     & -     & TL    &       & -     & -     & TL \\
    63    & tc1c80s8cf2-p16 &       & 32.530 & 50.0\% & TL    &       & \multicolumn{1}{r}{29.825} & \multicolumn{1}{r}{22.4\%} & TL    &       & \multicolumn{1}{r}{35.542} & \multicolumn{1}{r}{37.5\%} & TL    &       & -     & -     & TL \\
    64    & tc2c80s8cf3-p16 &       & 34.610 & 54.8\% & TL    &       & -     & -     & TL    &       & \multicolumn{1}{r}{33.142} & \multicolumn{1}{r}{36.0\%} & TL    &       & -     & -     & TL \\
    65    & tc2c80s8cf4-p16 &       & \multicolumn{1}{c}{-} & \multicolumn{1}{c}{-} & TL    &       & -     & -     & TL    &       & -     & -     & TL    &       & -     & -     & TL \\
    66    & tc0c80s8ct0-p16 &       & 38.075 & 45.3\% & TL    &       & -     & -     & TL    &       & \multicolumn{1}{r}{40.113} & \multicolumn{1}{r}{36.3\%} & TL    &       & -     & -     & TL \\
    67    & tc0c80s8ct1-p16 &       & 55.613 & 60.8\% & TL    &       & -     & -     & TL    &       & -     & -     & TL    &       & -     & -     & TL \\
    68    & tc1c80s8ct2-p16 &       & 34.580 & 53.3\% & TL    &       & \multicolumn{1}{r}{30.430} & \multicolumn{1}{r}{20.7\%} & TL    &       & \multicolumn{1}{r}{40.114} & \multicolumn{1}{r}{45.7\%} & TL    &       & -     & -     & TL \\
    69    & tc2c80s8ct3-p16 &       & 33.149 & 52.8\% & TL    &       & -     & -     & TL    &       & \multicolumn{1}{r}{32.911} & \multicolumn{1}{r}{36.2\%} & TL    &       & -     & -     & TL \\
    70    & tc2c80s8ct4-p16 &       & \multicolumn{1}{c}{-} & \multicolumn{1}{c}{-} & TL    &       & -     & -     & TL    &       & \multicolumn{1}{r}{63.465} & \multicolumn{1}{r}{59.7\%} & TL    &       & -     & -     & TL \\
    71    & tc0c80s12cf0-p24 &       & 37.099 & 43.7\% & TL    &       & -     & -     & TL    &       & \multicolumn{1}{r}{60.219} & \multicolumn{1}{r}{62.0\%} & TL    &       & -     & -     & TL \\
    72    & tc0c80s12cf1-p24 &       & 52.060 & 58.3\% & TL    &       & -     & -     & TL    &       & -     & -     & TL    &       & -     & -     & TL \\
    73    & tc1c80s12cf2-p24 &       & 33.727 & 51.6\% & TL    &       & \multicolumn{1}{r}{32.637} & \multicolumn{1}{r}{30.6\%} & TL    &       & \multicolumn{1}{r}{34.393} & \multicolumn{1}{r}{41.6\%} & TL    &       & \multicolumn{1}{r}{173.926} & \multicolumn{1}{r}{100.0\%} & TL \\
    74    & tc2c80s12cf3-p24 &       & 29.927 & 47.8\% & TL    &       & \multicolumn{1}{r}{29.803} & \multicolumn{1}{r}{34.8\%} & TL    &       & \multicolumn{1}{r}{37.849} & \multicolumn{1}{r}{46.5\%} & TL    &       & \multicolumn{1}{r}{195.129} & \multicolumn{1}{r}{100.0\%} & TL \\
    75    & tc2c80s12cf4-p24 &       & \multicolumn{1}{c}{-} & \multicolumn{1}{c}{-} & TL    &       & -     & -     & TL    &       & -     & -     & TL    &       & -     & -     & TL \\
    76    & tc0c80s12ct0-p24 &       & 36.673 & 43.3\% & TL    &       & \multicolumn{1}{r}{38.498} & \multicolumn{1}{r}{34.2\%} & TL    &       & \multicolumn{1}{r}{42.484} & \multicolumn{1}{r}{43.8\%} & TL    &       & -     & -     & TL \\
    77    & tc0c80s12ct1-p24 &       & 49.544 & 56.0\% & TL    &       & -     & -     & TL    &       & \multicolumn{1}{r}{53.002} & \multicolumn{1}{r}{49.3\%} & TL    &       & -     & -     & TL \\
    78    & tc1c80s12ct2-p24 &       & 31.347 & 47.8\% & TL    &       & \multicolumn{1}{r}{33.659} & \multicolumn{1}{r}{34.5\%} & TL    &       & \multicolumn{1}{r}{36.473} & \multicolumn{1}{r}{43.1\%} & TL    &       & \multicolumn{1}{r}{185.136} & \multicolumn{1}{r}{100.0\%} & TL \\
    79    & tc2c80s12ct3-p24 &       & 30.067 & 48.0\% & TL    &       & \multicolumn{1}{r}{29.838} & \multicolumn{1}{r}{26.8\%} & TL    &       & \multicolumn{1}{r}{36.161} & \multicolumn{1}{r}{49.6\%} & TL    &       & \multicolumn{1}{r}{192.812} & \multicolumn{1}{r}{100.0\%} & TL \\
    80    & tc2c80s12ct4-p24 &       & 54.834 & 72.3\% & TL    &       & -     & -     & TL    &       & -     & -     & TL    &       & -     & -     & TL \\
    \bottomrule
    \end{tabular}%
  \label{tab:addlabel}%
\end{table}%

\clearpage
\section{Results for the Linear Relaxation of each formulation}\label{App:C}

This section presents the results for the Linear Relaxation for each of the 80 instances and the four different models.

\begin{table}[htbp]
  \centering
  \caption{Linear Relaxation for the 10 customer instances}
  \scriptsize
    \begin{tabular}{cccccc}
    \toprule
    \textbf{\#} & \textbf{Instance} & \textbf{M1} & \textbf{M2} & \textbf{M3} & \textbf{M4} \\
    \midrule
    1     & tc0c10s2cf1-p4 & 5.808 & 7.506 & 7.291 & 7.348 \\
    2     & tc0c10s2ct1-p4 & 5.057 & 6.894 & 7.322 & 7.310 \\
    3     & tc0c10s3cf1-p6 & 5.223 & 6.467 & 7.205 & 7.205 \\
    4     & tc0c10s3ct1-p6 & 4.381 & 6.235 & 7.218 & 7.218 \\
    5     & tc1c10s2cf2-p4 & 3.864 & 4.765 & 5.571 & 5.571 \\
    6     & tc1c10s2cf3-p4 & 2.939 & 6.982 & 6.628 & 6.664 \\
    7     & tc1c10s2cf4-p4 & 5.253 & 6.928 & 7.226 & 7.226 \\
    8     & tc1c10s2ct2-p4 & 3.435 & 4.652 & 5.571 & 5.571 \\
    9     & tc1c10s2ct3-p4 & 2.743 & 5.896 & 6.481 & 6.481 \\
    10    & tc1c10s2ct4-p4 & 5.192 & 6.722 & 7.224 & 7.224 \\
    11    & tc1c10s3cf2-p6 & 4.230 & 5.293 & 5.571 & 5.571 \\
    12    & tc1c10s3cf3-p6 & 2.452 & 5.741 & 6.480 & 6.480 \\
    13    & tc1c10s3cf4-p6 & 4.917 & 6.727 & 7.348 & 7.342 \\
    14    & tc1c10s3ct2-p6 & 3.992 & 5.022 & 5.571 & 5.571 \\
    15    & tc1c10s3ct3-p6 & 2.301 & 5.727 & 6.480 & 6.480 \\
    16    & tc1c10s3ct4-p6 & 5.290 & 6.688 & 7.355 & 7.319 \\
    17    & tc2c10s2cf0-p4 & 1.034 & 3.728 & 3.797 & 3.797 \\
    18    & tc2c10s2ct0-p4 & 0.886 & 3.717 & 3.797 & 3.797 \\
    19    & tc2c10s3cf0-p6 & 1.020 & 3.709 & 3.797 & 3.797 \\
    20    & tc2c10s3ct0-p6 & 0.873 & 3.689 & 3.797 & 3.797 \\
    \midrule
          & Avg LR & 3.544 & 5.654 & 6.087 & 6.089 \\
    \bottomrule
    \end{tabular}%
  \label{tab:addlabel}%
\end{table}%

\begin{table}[htbp]
  \centering
    \caption{Linear Relaxation for the 20 customer instances}
  \scriptsize
    \begin{tabular}{rccccc}
    \toprule
    \multicolumn{1}{c}{\textbf{\#}} & \textbf{Instance} & \textbf{M1} & \textbf{M2} & \textbf{M3} & \textbf{M4} \\
    \midrule
    \multicolumn{1}{c}{21} & tc2c20s3cf0-p6 & 1.758 & 8.883 & 8.954 & 8.954 \\
    \multicolumn{1}{c}{22} & tc1c20s3cf1-p6 & 6.635 & 9.299 & 9.647 & 9.647 \\
    \multicolumn{1}{c}{23} & tc0c20s3cf2-p6 & 6.007 & 9.443 & 10.180 & 10.180 \\
    \multicolumn{1}{c}{24} & tc1c20s3cf3-p6 & 4.954 & 7.930 & 8.105 & 8.105 \\
    \multicolumn{1}{c}{25} & tc1c20s3cf4-p6 & 8.501 & 10.891 & 12.129 & 12.125 \\
    \multicolumn{1}{c}{26} & tc2c20s3ct0-p6 & 1.646 & 8.692 & 8.954 & 8.954 \\
    \multicolumn{1}{c}{27} & tc1c20s3ct1-p6 & 6.616 & 9.187 & 9.647 & 9.647 \\
    \multicolumn{1}{c}{28} & tc0c20s3ct2-p6 & 6.100 & 9.348 & 10.169 & 10.169 \\
    \multicolumn{1}{c}{29} & tc1c20s3ct3-p6 & 5.207 & 7.900 & 8.105 & 8.105 \\
    \multicolumn{1}{c}{30} & tc1c20s3ct4-p6 & 7.968 & 9.792 & 12.144 & 12.110 \\
    \multicolumn{1}{c}{31} & tc2c20s4cf0-p8 & 1.612 & 8.664 & 8.963 & 8.963 \\
    \multicolumn{1}{c}{32} & tc1c20s4cf1-p8 & 6.718 & 9.335 & 9.665 & 9.665 \\
    \multicolumn{1}{c}{33} & tc0c20s4cf2-p8 & 6.471 & 9.291 & 10.176 & 10.173 \\
    \multicolumn{1}{c}{34} & tc1c20s4cf3-p8 & 4.078 & 7.298 & 8.111 & 8.111 \\
    \multicolumn{1}{c}{35} & tc1c20s4cf4-p8 & 8.722 & 11.408 & 12.220 & 12.213 \\
    \multicolumn{1}{c}{36} & tc2c20s4ct0-p8 & 1.533 & 8.533 & 8.958 & 8.958 \\
    \multicolumn{1}{c}{37} & tc1c20s4ct1-p8 & 6.676 & 9.443 & 9.664 & 9.664 \\
    \multicolumn{1}{c}{38} & tc0c20s4ct2-p8 & 6.265 & 8.929 & 10.023 & 10.023 \\
    \multicolumn{1}{c}{39} & tc1c20s4ct3-p8 & 4.121 & 7.271 & 8.105 & 8.105 \\
    \multicolumn{1}{c}{40} & tc1c20s4ct4-p8 & 8.814 & 11.126 & 12.119 & 12.119 \\
    \midrule
          & Avg LR & 5.520 & 9.133 & 9.802 & 9.799 \\
    \bottomrule
    \end{tabular}%
  \label{tab:addlabel}%
\end{table}%

\begin{table}[htbp]
  \centering
      \caption{Linear Relaxation for the 40 customer instances}
  \scriptsize
    \begin{tabular}{rccccc}
    \toprule
    \multicolumn{1}{c}{\textbf{\#}} & \textbf{Instance} & \textbf{M1} & \textbf{M2} & \textbf{M3} & \textbf{M4} \\
    \midrule
    \multicolumn{1}{c}{41} & tc2c40s5cf2-p10 & 4.367 & 12.190 & 13.432 & 13.432 \\
    \multicolumn{1}{c}{42} & tc2c40s5cf3-p10 & 2.917 & 10.687 & 10.801 & 10.835 \\
    \multicolumn{1}{c}{43} & tc0c40s5cf4-p10 & 10.779 & 15.338 & 16.409 & 16.409 \\
    \multicolumn{1}{c}{44} & tc2c40s8cf2-p16 & 3.992 & 11.742 & 13.420 & 13.420 \\
    \multicolumn{1}{c}{45} & tc2c40s8cf3-p16 & 2.794 & 10.813 & 10.881 & 10.928 \\
    \multicolumn{1}{c}{46} & tc0c40s8ct0-p16 & 9.831 & 15.089 & 17.484 & 17.484 \\
    \multicolumn{1}{c}{47} & tc2c40s8ct3-p16 & 2.773 & 10.686 & 10.835 & 10.851 \\
    \multicolumn{1}{c}{48} & tc0c40s5cf0-p10 & 11.549 & 17.638 & 17.520 & 17.524 \\
    \multicolumn{1}{c}{49} & tc1c40s5cf1-p10 & 8.751 & 16.654 & 18.917 & 18.917 \\
    \multicolumn{1}{c}{50} & tc0c40s5ct0-p10 & 11.518 & 17.273 & 17.496 & 17.496 \\
    \multicolumn{1}{c}{51} & tc1c40s5ct1-p10 & 9.854 & 18.006 & 18.931 & 18.931 \\
    \multicolumn{1}{c}{52} & tc2c40s5ct2-p10 & 4.205 & 13.042 & 13.426 & 13.426 \\
    \multicolumn{1}{c}{53} & tc2c40s5ct3-p10 & 2.917 & 10.337 & 10.756 & 10.756 \\
    \multicolumn{1}{c}{54} & tc0c40s5ct4-p10 & 10.834 & 15.103 & 16.341 & 16.342 \\
    \multicolumn{1}{c}{55} & tc0c40s8cf0-p16 & 10.660 & 15.656 & 17.488 & 17.488 \\
    \multicolumn{1}{c}{56} & tc1c40s8cf1-p16 & 9.076 & 16.537 & 18.917 & 18.917 \\
    \multicolumn{1}{c}{57} & tc0c40s8cf4-p16 & 10.216 & 14.298 & 16.205 & 16.205 \\
    \multicolumn{1}{c}{58} & tc1c40s8ct1-p16 & 9.512 & 16.954 & 18.917 & 18.917 \\
    \multicolumn{1}{c}{59} & tc2c40s8ct2-p16 & 3.725 & 12.444 & 13.420 & 13.420 \\
    \multicolumn{1}{c}{60} & tc0c40s8ct4-p16 & 10.648 & 14.621 & 16.198 & 16.198 \\
    \midrule
          & Avg LR & 7.546 & 14.255 & 15.390 & 15.395 \\
    \bottomrule
    \end{tabular}%
  \label{tab:addlabel}%
\end{table}%

\begin{table}[htbp]
  \centering
\caption{Linear Relaxation for the 80 customer instances}
  \scriptsize
    \begin{tabular}{cccccc}
    \toprule
    \textbf{\#} & \textbf{Instance} & \textbf{M1} & \textbf{M2} & \textbf{M3} & \textbf{M4} \\
    \midrule
    61    & tc0c80s8cf0-p16 & 15.852 & 22.286 & 22.929 & 22.931 \\
    62    & tc0c80s8cf1-p16 & 17.589 & 26.468 & 26.861 & 26.861 \\
    63    & tc1c80s8cf2-p16 & 10.981 & 19.097 & 19.710 & 19.716 \\
    64    & tc2c80s8cf3-p16 & 7.123 & 17.025 & 17.425 & 17.425 \\
    65    & tc2c80s8cf4-p16 & 8.220 & 21.759 & 22.744 & 22.744 \\
    66    & tc0c80s8ct0-p16 & 15.927 & 22.251 & 22.921 & 22.924 \\
    67    & tc0c80s8ct1-p16 & 17.482 & 26.445 & 26.872 & 26.872 \\
    68    & tc1c80s8ct2-p16 & 10.956 & 19.034 & 19.647 & 19.647 \\
    69    & tc2c80s8ct3-p16 & 7.175 & 17.133 & 17.425 & 17.425 \\
    70    & tc2c80s8ct4-p16 & 7.799 & 21.930 & 22.744 & 22.744 \\
    71    & tc0c80s12cf0-p24 & 14.626 & 21.332 & 22.884 & 22.884 \\
    72    & tc0c80s12cf1-p24 & 17.005 & 25.285 & 26.847 & 26.847 \\
    73    & tc1c80s12cf2-p24 & 10.695 & 18.832 & 19.669 & 19.669 \\
    74    & tc2c80s12cf3-p24 & 6.930 & 16.698 & 17.424 & 17.424 \\
    75    & tc2c80s12cf4-p24 & 8.019 & 21.344 & 22.744 & 22.744 \\
    76    & tc0c80s12ct0-p24 & 14.736 & 21.174 & 22.884 & 22.884 \\
    77    & tc0c80s12ct1-p24 & 17.047 & 25.462 & 26.847 & 26.847 \\
    78    & tc1c80s12ct2-p24 & 10.941 & 19.012 & 19.647 & 19.647 \\
    79    & tc2c80s12ct3-p24 & 6.869 & 16.738 & 17.424 & 17.424 \\
    80    & tc2c80s12ct4-p24 & 7.691 & 21.666 & 22.744 & 22.744 \\
    \midrule
          & Avg LR & 11.683 & 21.049 & 21.920 & 21.920 \\
    \bottomrule
    \end{tabular}%
  \label{tab:addlabel}%
\end{table}%

\clearpage
\section{Results considering the best objective function and the best lower bound}\label{App:D}

This section presents the results for the 80 instances considering the best objective function and the best lower bound obtained by the Gurobi solver in the four different formulations.

\begin{table}[htbp]
  \centering
  \caption{10 customer results considering the best objective function and the best lower bound}
    \scriptsize
    \begin{tabular}{ccccc}
    \toprule
    \textbf{\#} & \textbf{Instance} & \textbf{Obj Min} & \textbf{LB Max} & \textbf{Gap} \\
    \midrule
    1     & tc0c10s2cf1-p4 & 19.753 & 19.753 & 0.0\% \\
    2     & tc0c10s2ct1-p4 & 11.378 & 11.378 & 0.0\% \\
    3     & tc0c10s3cf1-p6 & 10.921 & 10.921 & 0.0\% \\
    4     & tc0c10s3ct1-p6 & 10.536 & 10.536 & 0.0\% \\
    5     & tc1c10s2cf2-p4 & 9.034 & 9.034 & 0.0\% \\
    6     & tc1c10s2cf3-p4 & 12.583 & 12.583 & 0.0\% \\
    7     & tc1c10s2cf4-p4 & 16.097 & 16.097 & 0.0\% \\
    8     & tc1c10s2ct2-p4 & 9.199 & 9.199 & 0.0\% \\
    9     & tc1c10s2ct3-p4 & 12.307 & 12.307 & 0.0\% \\
    10    & tc1c10s2ct4-p4 & 13.826 & 13.826 & 0.0\% \\
    11    & tc1c10s3cf2-p6 & 9.034 & 9.034 & 0.0\% \\
    12    & tc1c10s3cf3-p6 & 12.583 & 12.583 & 0.0\% \\
    13    & tc1c10s3cf4-p6 & 14.902 & 14.902 & 0.0\% \\
    14    & tc1c10s3ct2-p6 & 9.199 & 9.199 & 0.0\% \\
    15    & tc1c10s3ct3-p6 & 12.932 & 12.932 & 0.0\% \\
    16    & tc1c10s3ct4-p6 & 13.205 & 13.205 & 0.0\% \\
    17    & tc2c10s2cf0-p4 & 11.345 & 11.345 & 0.0\% \\
    18    & tc2c10s2ct0-p4 & 11.345 & 11.345 & 0.0\% \\
    19    & tc2c10s3cf0-p6 & 11.285 & 11.285 & 0.0\% \\
    20    & tc2c10s3ct0-p6 & 11.285 & 11.285 & 0.0\% \\
    \bottomrule
    \end{tabular}%
  \label{tab:addlabel}%
\end{table}%

\begin{table}[htbp]
  \centering
    \caption{20 customer results considering the best objective function and the best lower bound}
    \scriptsize
    \begin{tabular}{ccccc}
    \toprule
    \textbf{\#} & \textbf{Instance} & \textbf{Obj Min} & \textbf{LB Max} & \textbf{Gap} \\
    \midrule
    21    & tc2c20s3cf0-p6 & 24.040 & 16.385 & 31.8\% \\
    22    & tc1c20s3cf1-p6 & 14.747 & 13.799 & 6.4\% \\
    23    & tc0c20s3cf2-p6 & 16.651 & 14.758 & 11.4\% \\
    24    & tc1c20s3cf3-p6 & 11.861 & 11.861 & 0.0\% \\
    25    & tc1c20s3cf4-p6 & 16.088 & 16.087 & 0.0\% \\
    26    & tc2c20s3ct0-p6 & 24.275 & 21.804 & 10.2\% \\
    27    & tc1c20s3ct1-p6 & 15.083 & 13.783 & 8.6\% \\
    28    & tc0c20s3ct2-p6 & 16.627 & 14.669 & 11.8\% \\
    29    & tc1c20s3ct3-p6 & 11.861 & 11.861 & 0.0\% \\
    30    & tc1c20s3ct4-p6 & 15.792 & 15.791 & 0.0\% \\
    31    & tc2c20s4cf0-p8 & 24.672 & 18.373 & 25.5\% \\
    32    & tc1c20s4cf1-p8 & 16.041 & 13.410 & 16.4\% \\
    33    & tc0c20s4cf2-p8 & 16.206 & 14.353 & 11.4\% \\
    34    & tc1c20s4cf3-p8 & 12.967 & 12.967 & 0.0\% \\
    35    & tc1c20s4cf4-p8 & 15.825 & 15.824 & 0.0\% \\
    36    & tc2c20s4ct0-p8 & 26.018 & 18.361 & 29.4\% \\
    37    & tc1c20s4ct1-p8 & 16.070 & 14.925 & 7.1\% \\
    38    & tc0c20s4ct2-p8 & 16.056 & 14.176 & 11.7\% \\
    39    & tc1c20s4ct3-p8 & 12.866 & 12.866 & 0.0\% \\
    40    & tc1c20s4ct4-p8 & 15.825 & 15.824 & 0.0\% \\
    \bottomrule
    \end{tabular}%
  \label{tab:addlabel}%
\end{table}%

\begin{table}[htbp]
  \centering
      \caption{40 customer results considering the best objective function and the best lower bound}
    \scriptsize
    \begin{tabular}{ccccc}
    \toprule
    \textbf{\#} & \textbf{Instance} & \textbf{Obj Min} & \textbf{LB Max} & \textbf{Gap} \\
    \midrule
    41    & tc2c40s5cf2-p10 & 27.425 & 20.268 & 26.1\% \\
    42    & tc2c40s5cf3-p10 & 19.171 & 15.007 & 21.7\% \\
    43    & tc0c40s5cf4-p10 & 33.034 & 21.209 & 35.8\% \\
    44    & tc2c40s8cf2-p16 & 26.705 & 20.035 & 25.2\% \\
    45    & tc2c40s8cf3-p16 & 19.642 & 15.251 & 22.4\% \\
    46    & tc0c40s8ct0-p16 & 24.802 & 22.268 & 10.2\% \\
    47    & tc2c40s8ct3-p16 & 22.545 & 15.017 & 33.4\% \\
    48    & tc0c40s5cf0-p10 & 27.245 & 23.216 & 14.8\% \\
    49    & tc1c40s5cf1-p10 & -     & 29.004 & - \\
    50    & tc0c40s5ct0-p10 & 26.186 & 22.735 & 13.2\% \\
    51    & tc1c40s5ct1-p10 & - & 28.613 & - \\
    52    & tc2c40s5ct2-p10 & 26.216 & 19.894 & 24.1\% \\
    53    & tc2c40s5ct3-p10 & 18.714 & 14.716 & 21.4\% \\
    54    & tc0c40s5ct4-p10 & 31.707 & 20.929 & 34.0\% \\
    55    & tc0c40s8cf0-p16 & 25.930 & 23.298 & 10.2\% \\
    56    & tc1c40s8cf1-p16 & -     & 26.090 & - \\
    57    & tc0c40s8cf4-p16 & 28.227 & 20.846 & 26.1\% \\
    58    & tc1c40s8ct1-p16 & -     & 26.309 & - \\
    59    & tc2c40s8ct2-p16 & 26.107 & 19.519 & 25.2\% \\
    60    & tc0c40s8ct4-p16 & 28.502 & 20.888 & 26.7\% \\
    \bottomrule
    \end{tabular}%
  \label{tab:addlabel}%
\end{table}%

\begin{table}[htbp]
  \centering
    \caption{80 customer results considering the best objective function and the best lower bound}
    \scriptsize
    \begin{tabular}{ccccc}
    \toprule
    \textbf{\#} & \textbf{Instance} & \textbf{Obj Min} & \textbf{LB Max} & \textbf{Gap} \\
    \midrule
    61    & tc0c80s8cf0-p16 & 38.205 & 28.873 & 24.4\% \\
    62    & tc0c80s8cf1-p16 & 55.182 & 32.653 & 40.8\% \\
    63    & tc1c80s8cf2-p16 & 29.825 & 24.519 & 17.8\% \\
    64    & tc2c80s8cf3-p16 & 34.610 & 23.230 & 32.9\% \\
    65    & tc2c80s8cf4-p16 & -     & 30.144 & - \\
    66    & tc0c80s8ct0-p16 & 38.075 & 28.619 & 24.8\% \\
    67    & tc0c80s8ct1-p16 & 55.613 & 32.638 & 41.3\% \\
    68    & tc1c80s8ct2-p16 & 30.430 & 24.119 & 20.7\% \\
    69    & tc2c80s8ct3-p16 & 32.911 & 21.568 & 34.5\% \\
    70    & tc2c80s8ct4-p16 & 63.465 & 29.096 & 54.2\% \\
    71    & tc0c80s12cf0-p24 & 37.099 & 28.423 & 23.4\% \\
    72    & tc0c80s12cf1-p24 & 52.060 & 31.746 & 39.0\% \\
    73    & tc1c80s12cf2-p24 & 32.637 & 23.848 & 26.9\% \\
    74    & tc2c80s12cf3-p24 & 29.803 & 22.296 & 25.2\% \\
    75    & tc2c80s12cf4-p24 & -     & 29.045 & - \\
    76    & tc0c80s12ct0-p24 & 36.673 & 28.072 & 23.5\% \\
    77    & tc0c80s12ct1-p24 & 49.544 & 31.571 & 36.3\% \\
    78    & tc1c80s12ct2-p24 & 31.347 & 23.436 & 25.2\% \\
    79    & tc2c80s12ct3-p24 & 29.838 & 21.845 & 26.8\% \\
    80    & tc2c80s12ct4-p24 & 54.834 & 28.629 & 47.8\% \\
    \bottomrule
    \end{tabular}%
  \label{tab:addlabel}%
\end{table}%

\end{document}